\documentstyle[11pt,amssymb,amsfonts]{article}

\begin{document}
\newcommand{\p}{\parallel }
\makeatletter \makeatother
\newtheorem{th}{Theorem}[section]
\newtheorem{lem}{Lemma}[section]
\newtheorem{de}{Definition}[section]
\newtheorem{rem}{Remark}[section]
\newtheorem{cor}{Corollary}[section]
\renewcommand{\theequation}{\thesection.\arabic {equation}}

\title{{\bf The Noncommutative Infinitesimal Equivariant Index Formula}}

\author{ Yong Wang }

\date{}
\maketitle

\begin{abstract}
In this paper, we establish an infinitesimal equivariant index formula in the
noncommutative geometry framework using Greiner's approach to heat kernel
 asymptotics. An infinitesimal equivariant index formula for odd dimensional manifolds is also given. We define infinitesimal equivariant eta cochains, prove their regularity and give an explicit formula for them.
We also establish an infinitesimal equivariant family index formula and introduce the infinitesimal equivariant eta forms as well as compare them with the equivariant eta forms.
\\

\noindent{\bf Keywords:}\quad Infinitesimal equivariant Chern-Connes
characters; infinitesimal equivariant eta cochains; infinitesimal
equivariant family index;
 infinitesimal equivariant eta forms
\\

\noindent {\bf MSC(2010)}: 58J20, 19K56\\

\end{abstract}

\section{Introduction}
    \quad
 The Atiyah-Bott-Segal-Singer index formula is a
generalization of the Atiyah-Singer index theorem to manifolds admitting group actions.
In [BV1], Berline and Vergne gave a heat kernel proof of the
Atiyah-Bott-Segal-Singer index formula. In [LYZ], Lafferty, Yu and
Zhang gave a very simple and direct geometric proof to the
equivariant index formula of Dirac operators. In [PW], Ponge and Wang
gave another proof of the equivariant index formula using
Greiner's approach to heat kernel asymptotics. For manifolds with boundary,  Donnelly [Do] introduced the equivariant eta
invariant and generalized the Atiyah-Patodi-Singer index theorem to the equivariant setting. Zhang proved the regularity of the equivariant eta invariant in [Zh].  In [Fa], Fang established an equivariant index formula for odd dimensional manifolds.\\
 \indent The equivariant index formula has an infinitesimal version,
which is called the Kirillov formula. Berline and Vergne
[BV2] established the Kirillov formula using the equivariant index formula and the localization
formula. Bismut introduced the Bismut Laplacian and gave a direct heat kernel proof of the Kirillov
formula in [Bi]. The infinitesimal equivariant index formula for manifolds with boundary
was established in [Go], where Goette introduced infinitesimal
equivariant eta invariants and compared equivariant eta invariants with infinitesimal equivariant eta invariants.\\
\indent On the noncommutative geometry side, Connes [Co] defined the Chern-Connes
character of a $\theta$-summable Fredholm module $(H,D)$ over a
unital $C^*$-algebra $A$, which takes value in the entire cyclic
cohomology of $A$. In [JLO], Jaffe, Lesniewski
and Osterwalder introduced an equivariant but convenient version of the Chern-Connes
character, which is known as the JLO character. The JLO character was computed in [CM1] and [BlF].
An explicit formula of the equivariant JLO character associated to the
invariant Dirac operator, in the presence of a countable discrete
group action on a smooth compact spin Riemannian manifold, was given
by Azmi [Az] and moreover it was shown that this equivariant
cocycle is an element of the delocalized cohomology, paired
with an equivariant K-theory idempotent. When $G$ is a
compact Lie group, Chern and Hu [CH] gave an explicit formula of the
equivariant Chern-Connes character associated to a G-equivariant
$\theta$-summable Fredholm module. In [Ge1], for odd dimensional
manifolds, the spectral flow was written as pairing of the JLO
character with the odd Chern character of an idempotent matrix.

In the framework of noncommutative geometry, Wu established an Atiyah-Patodi-Singer index theorem in [Wu]. To do so, he introduced the total eta
invariant (called the higher eta invariant in [Wu]), which is a
generalization of the classical Atiyah-Patodi-Singer eta invariants.
Wu then proved its regularity using the Getzler symbol calculus as
adopted in [BlF] and computed its radius of convergence. Subsequently, he established a variation
formula of eta cochains, which he used to obtain the noncommutative
Atiyah-Patodi-Singer index theorem. In [Ge2], using superconnection,
Getzler gave another proof to the noncommutative
Atiyah-Patodi-Singer index theorem, which is more difficult but
avoided mention of the operators $b$ and $B$ in cyclic cohomology.
In [Wa1], we introduced equivariant eta chains and established an
equivariant noncommutative Atiyah-Patodi-Singer index formula which generalized
Wu's theorem to the equivariant setting.

This paper is devoted to establish an infinitesimal
equivariant index formula in the noncommutative geometry framework
using Greiner's approach to heat kernel asymptotics as well as establish an
infinitesimal equivariant index formula for odd dimensional manifolds. In the same framework, we also give an infinitesimal equivariant index formula for
manifolds with boundary.\\
 \indent Let $D$ be a differential operator acting on a fiber bundle $M$ over a compact
space $B.$ If $D$ is elliptic along the fibers, then $D$ can be viewed as a family of elliptic operators parameterized by $B.$
 Atiyah and Singer defined a more general index for $D$ which
is an element in the $K$ group $K(B).$  This index is called family index.
Atiyah and Singer proved that the analytic and topological indices coincide in $K(B).$ As a consequence,
 they could determine the Chern character of the difference
bundle ${\rm Ker}D-{\rm Coker}D$ and gave a cohomology expression of the Chern character of the difference
bundle in terms of certain characteristic classes using Chern-Weil's theory. A nice exposition of  family
index theory can be found in [BGV, Ch.10]. In order
to prove family rigidity theorems for certain elliptic operators, Liu and Ma established an equivariant
family index formula [LM]. In [Wa2], using Greiner's approach to
heat kernel asymptotics, we gave another proof of the local equivariant index
theorem for a family of Dirac operators. We also introduced the equivariant eta
forms and proved their regularity in [Wa2]. The current paper will study the infinitesimal versions too, including an infinitesimal equivariant family index formula, the definition of infinitesimal equivariant eta forms and the comparison of
them with equivariant eta forms.

This paper is organized as follows: In Section 2, we establish an infinitesimal equivariant index formula in the
noncommutative geometry framework using Greiner's approach to heat
kernel asymptotics. An infinitesimal equivariant index formula
for odd dimensional manifolds is also established. In Section 3, we define
truncated infinitesimal equivariant eta cochains and prove their
regularity as well as give a formula for them. In Section 4, a proof of an infinitesimal equivariant
family index formula is given. We also introduce infinitesimal equivariant
eta forms and compare them with equivariant eta forms.\\

\section{ The noncommutative infinitesimal equivariant index formula
 }

 {\bf 2.1 The infinitesimal equivariant JLO cocycle}\\

 \indent Let $M$ be a compact
oriented even dimensional Riemannian manifold without boundary with a
fixed spin structure and $S$ be the bundle of spinors on $M$. Denote
by $D$ the associated Dirac operator on $H=L^2(M;S)$, the Hilbert
space of $L^2$-sections of the bundle $S$. Let $c(df): S\rightarrow
S$ denote the Clifford action with $f\in C^{\infty}(M)$. Suppose
that $G$ is a compact connected Lie group acting on $M$ by
orientation-preserving isometries preserving the spin structure and $\mathfrak{g}$ is the Lie algebra of $G$. Then $G$ commutes with the Dirac operator. For $X\in\mathfrak{g}$, let $X_M(p)=\frac{d}{dt}|_{t=0}e^{-tX}p$ be the Killing
field induced by $X$, Let $c(X)$ denote the Clifford action by $X_M$, and ${\mathfrak{L}}_X$ denote the Lie derivative. Define $\mathfrak{g}$-equivariant modifications of $D$ and $D^2$ for $X\in {\mathfrak{g}}$ as follows:\\
$$D_X:=D-\frac{1}{4}c(X);~~H_X:=D^2_{-X}+{\mathfrak{L}}_X=(D+\frac{1}{4}c(X))^2+{\mathfrak{L}}_X,\eqno(2.1)$$
then $H_X$ is the equivariant Bismut Laplacian. Let ${\mathbb{C}}[{\mathfrak{g}}^*]$ denote the space of formal power series in $X\in {\mathfrak{g}}$ and $\psi_t$ be the rescaling operator on ${\mathbb{C}}[{\mathfrak{g}}^*]$ defined by $X\rightarrow \frac{X}{t}$ for $t>0$.\\
 \indent Let
$$A=C_G^{\infty}(M)=\{f\in C^{\infty}(M)|f(g\cdot x)=f(x), g\in~ G, x\in~ M\},$$
 then the data $(A,H,D+\frac{1}{4}c(X),G)$ defines a non selfadjoint perturbation of finitely summable
(hence $\theta$-summable) equivariant unbounded Fredholm module
$(A,H,D,G)$ in the sense of [KL] (for details, see [CH] and [KL]).
For $(A,H,D+\frac{1}{4}c(X),G)$,
  {\bf The truncated infinitesimal equivariant JLO cochain} ${{\bf
{\rm ch}}}_{
2k}(\sqrt{t}D,X)_J$ can be defined by the formula:\\

$${{\bf {\rm ch}}_{2k}}(\sqrt{t}D,X)(f^0,\cdots ,
f^{2k})_J:
=t^k\int_{\triangle_{2k}}{\rm
Str}\left[\psi_t e^{-t{
\mathfrak{L}}_X}f^0e^{-\sigma_0t(D+\frac{1}{4}c(X))^2}c(df^1)\right.$$
$$\left.\cdot e^{-\sigma_1t(D+\frac{1}{4}c(X))^2}\cdots c(df^{2k})e^{-\sigma_{2k}t(D+\frac{1}{4}c(X))^2}\right]_Jd{\rm Vol}_{\Delta_{2k}},\eqno(2.2)$$
 where
$\triangle_{2k}=\{(\sigma_0,\cdots,\sigma_{2k})|~\sigma_0+\cdots \sigma_{2k}=1\}$ is the $2k$-simplex.
For an integer $J\geq 0$, denote by ${\mathbb{C}}[{\mathfrak{g}}^*]_J$ the space of polynomials in $X\in {\mathfrak{g}}$ of degree
$\leq J$ and let $(\cdot)_J:~ {\mathbb{C}}[{\mathfrak{g}}^*]\rightarrow {\mathbb{C}}[{\mathfrak{g}}^*]_J$ be the natural projection. Then
${{\bf {\rm ch}}_{2k}}(\sqrt{t}D,X)(f^0,\cdots ,
f^{2k})_J$ is controlled by $\psi_t(C(X))_Jt^{k}{\rm Tr}(e^{-\frac{t}{2}D^2})$ via the following lemma 2.2 (similar to Lemma 2.1 in
 [GS]), so it is well-defined for $t\in (0,+\infty)$.
We will
compute the limit of the $J$-jet of the infinitesimal equivariant JLO cochain
$${\rm lim}_{t\rightarrow 0}{{\bf {\rm ch}}_{2k}}(\sqrt{t}D,X)(f^0,\cdots ,
f^{2k})_J.$$

\indent In the following, we give
some estimates about ${{\bf {\rm ch}}_{2k}}(\sqrt{t}D,X)(f^0,\cdots ,
f^{2k})_J.$ Let $H$ be a Hilbert space. For
$q\geq0$, denote by $||.||_q$ the Schatten $p$-norm on the Schatten
ideal
$L^p$. Let $L(H)$ denote the Banach algebra of bounded operators on $H$.\\

\noindent {\bf Lemma 2.1}~([CH],[Fe]){\it~~(i)~~${\rm Tr}(AB)={\rm
Tr}(BA)$, for $A,~B\in L(H)$ and $AB, ~BA\in
L^1$.\\
~~(ii)~~For $A\in L^1,$ we have $|{\rm Tr}(A)|\leq ||A||_1$,
$||A||\leq ||A||_1$.\\
~~(iii)~~For $A\in L^q$ and $B\in L(H)$, we have: $||AB||_q\leq
||B||||A||_q$, $||BA||_q\leq ||B||||A||_q$.\\
 ~~(iv)~(H\"{o}lder
Inequality)~~If $\frac{1}{r}=\frac{1}{p}+\frac{1}{q},~p,q,r>0,~A\in
L^p,~B\in
L^q,$ then $AB\in L^r$ and $||AB||_r\leq ||A||_p||B||_q$.}\\

Fix basis $e_1,\cdots,e_n$ of ${\mathfrak{g}}$ and let $X=x_1e_1+\cdots x_ne_n.$ A $J$-degree polynomial on $X$ means a $J$-degree polynomial on $x_1,\cdots,x_n.$\\

\noindent{\bf Lemma 2.2} {\it For any $1\geq u> 0$,~$t>0$ and $t$ is
small, $X\in {\mathfrak{g}}$ and any order $l$
differential operator ${B}$, we have:}\\
$$||e^{-utH_X}_J{B}||_{u^{-1}}\leq C(X)_Ju^{-\frac{l}{2}}t^{-\frac{l}{2}}({\rm
tr}[e^{-\frac{tD^2}{2}}])^u,\eqno(2.3)$$
{\it where $C(X)_J$ is a $J$-degree polynomial with constant coefficients on $X$.} \\

\noindent {\bf Proof.} Let $H_X=D^2+F_X$, where $F_X$ is a first order differential operator with degree $\geq 1$ coefficients depending on $X$. By the Duhamel principle, it is that
$$||e^{-utH_X}_J{B}||_{u^{-1}}=||
\sum^J_{m\geq
0}(-ut)^m\int_{\triangle_m}e^{-v_0utD^2}F_Xe^{-v_1utD^2}$$
$$\cdot F_X\cdots e^{-v_{m-1}utD^2}F_Xe^{-v_mutD^2}{B}dv||_{u^{-1}}.\eqno(2.4)$$
We estimate the term for $m=2$ in the right hand side of (2.4), and other
terms can be estimated similarly. We split $\triangle_2=J_0\cup J_1 \cup J_2$
where $J_i=\{(v_0,v_1,v_2)\in \Delta_2|v_i\geq \frac{1}{3}\}.$ Then,
\begin{eqnarray*}
&&(ut)^2||\int_{J_0}e^{-v_0utD^2}F_{X}e^{-v_1utD^2}
F_{X}e^{-v_2utD^2}{B}dv||_{u^{-1}}\\
&&\leq (ut)^2\int_{J_0}||e^{-\frac{v_0ut}{2}D^2}||_{(uv_0)^{-1}}
||e^{-\frac{v_0ut}{2}D^2}(1+D^2)^{\frac{l+2}{2}}||||(1+D^2)^{-\frac{l+2}{2}}
F_{X}(1+D^2)^{\frac{l+1}{2}}||\\
&&\cdot||e^{-{v_1ut}D^2}||_{(uv_1)^{-1}}||(1+D^2)^{-\frac{l+1}{2}}
F_{X}(1+D^2)^{\frac{l}{2}}||||e^{-{v_2ut}D^2}||_{(uv_2)^{-1}}||1+D^2)^{-\frac{l}{2}}{B}||dv\\
&&\leq (ut)^2\int_{J_0} \left({\rm
Tr}e^{-\frac{t}{2}D^2}\right)^{uv_0}\left({\rm
Tr}e^{-{t}D^2}\right)^{u(v_1+v_2)}(uv_0t)^{-\frac{l+2}{2}}\\
&&\cdot||(1+D^2)^{-\frac{l+2}{2}}
F_{X}(1+D^2)^{\frac{l+1}{2}}||||(1+D^2)^{-\frac{l+1}{2}}
F_{X}(1+D^2)^{\frac{l}{2}}||||1+D^2)^{-\frac{l}{2}}{B}||dv\\
&&\leq
 C(X)_2\left({\rm
Tr}e^{-\frac{t}{2}D^2}\right)^{u}(ut)^{-\frac{l}{2}+1},~~~~~~~~~~(2.5)
\end{eqnarray*}
where we use that $F_{X}$ is a first order differential
operator and the equality
$${\rm
sup}\{(1+x)^{\frac{l}{2}}e^{-\frac{utx}{2}}\}=(ut)^{-\frac{l}{2}}e^{-\frac{l-ut}{2}}.\eqno(2.6)$$
 $J_1$ and $J_2$ can be estimated similarly. For
the general $m$, we get
$$||(-ut)^m\int_{\triangle_m}e^{-v_0utD^2}F_{X}e^{-v_1utD^2}F_{X}\cdots
e^{-v_{m-1}utD^2}$$
$$\cdot F_{X}e^{-v_mutD^2}{B}dv||_{u^{-1}}\leq C_2\left({\rm
Tr}e^{-\frac{t}{2}D^2}\right)^{u}(ut)^{-\frac{l}{2}+\frac{m}{2}}.\eqno(2.7)$$
By (2.4) and (2.7), (2.3) is obtained. $\Box$\\

\indent Similarly to Lemmas 4.3 and 4.4 in [Wa2], we have\\

 \noindent{\bf Lemma 2.3}{\it ~Let ${B}_1,~{B}_2$ be positive order
 $p,~q$ pseudodifferential operators respectively, then for any
 $s,~t>0,~0\leq u\leq1$ and $t$ is small, $X\in {\mathfrak{g}}$, we have the following estimate:}\\
$$||[{B}_1e^{-ustH_X}{B}_2e^{-(1-u)stH_X}]_J||_{s^{-1}}\leq C(X)_Js^{-\frac{p+q}{2}}t^{-\frac{p+q}{2}}({\rm
tr}[e^{-\frac{tD^2}{4}}])^s.\eqno(2.8)$$\\\

\indent Let ${B}$ be an
operator  and $l$ be a positive interger.
 Write
$${B}^{[l]}=[H_X,{B}^{[l-1]}],~{B}^{[0]}={B}.$$\\

\noindent{\bf Lemma 2.4}~{\it Let ${B}$ be a finite order
differential
operator with coefficients on $X$ , then for any $s>0$, we have:}\\
$$[e^{-sH_X}{B}]_J=\sum^{N-1}_{l=0}\frac{(-1)^l}{l!}s^l[{B}^{[l]}e^{-sH_X}]_J+(-1)^Ns^N({B}^{[N]}(s))_J,\eqno(2.9)$$
\noindent {\it where ${B}^{[N]}(s)$ is given by}\\
$${B}^{[N]}(s)=\int_{\triangle_N}e^{-u_1sH_X}{B}^{[N]}e^{-(1-u_1)sH_X}du_1du_2\cdots du_N.\eqno(2.10)$$
\\

\noindent{\bf Lemma 2.5}~{\it Let ${B}$ be a finite order
differential
operator with coefficients on $X$ , then for any $s>0$, we have:}\\
$$[Be^{-sH_X}]_J=\sum^{N-1}_{l=0}\frac{(-1)^l}{l!}s^l[e^{-sH_X}{B}^{[l]}]_J+(-1)^Ns^N({B}_1^{[N]}(s))_J,\eqno(2.11)$$
\noindent {\it where ${B}_1^{[N]}(s)$ is given by}\\
$${B}_1^{[N]}(s)=\int_{\triangle_N}e^{-(1-u_1)sH_X}{B}^{[N]}e^{-u_1sH_X}du_1du_2\cdots du_N.\eqno(2.12)$$
\\

  \indent Since ${\mathfrak{L}}_X$ commutes with $D$, $c(X)$ and $f\in C_G^{\infty}(M)$, then
by Lemma 2.4, we have:\\
$$[e^{-t{\mathfrak{L}}_X}f^0e^{-s_1t(D+\frac{1}{4}c(X))^2}c(df^1)
e^{-(s_2-s_1)t(D+\frac{1}{4}c(X))^2}\cdots c(df^{2k})e^{-(1-s_{2k})t(D+\frac{1}{4}c(X))^2}]_J$$
$$=\sum^{N-1}_{\lambda _{1},\cdots ,\lambda _{2k}=0}\frac{(-1)^{\lambda _1+\cdots +\lambda _{2k}}{s_1}^{\lambda _{1}}\cdots
 s_{2k}^{\lambda
_{2k}}t^{\lambda _{1}+\cdots +\lambda _{2k}}} {\lambda
_{1}!\cdots\lambda _{2k}!}[ f^0[c(df^1)]^{[\lambda _{1}]}\cdots
[c(df^{2k})]^{[\lambda _{2k}]}e^{-tH_X}]_J$$
$$+\sum_{1\leq q\leq {2k}}\sum^{N-1}_{\lambda _{1},\cdots ,\lambda _{q-1}=0}
\frac{(-1)^{\lambda _1+\cdots +\lambda _{q-1}+N}s_1^{\lambda
_1}\cdots s_{q-1}^{\lambda _{q-1}}s_q^Nt^{\lambda _1+\cdots +\lambda
_{q-1}+N}}{\lambda _{1}!\cdots\lambda
_{q-1}!}[f^0[c(df^1)]^{[\lambda _{1}]}$$
$$\cdots[c(df^{q-1})]^{[\lambda _{q-1}]}\{ [c(df^q)]^{[N]}(s_qt)\}
e^{-(s_{q+1}-s_q)tH_X}\cdots
c(df^{2k})e^{-(1-s_{2k})tH_X}]_J.\eqno(2.13)$$
Since $f^0[c(df^1)]^{[\lambda _{1}]}
\cdots[c(df^{q-1})]^{[\lambda _{q-1}]}$ is a $\lambda _{1}+\cdots+\lambda _{q-1}$ order differential operator, we get by Lemma 2.2 and Lemma 2.3 (see pp. 61-62 in [Fe]) that
$$\left|\psi_t\int_{\Delta_{2k}}t^k
\sum_{1\leq q\leq {2k}}\sum^{N-1}_{\lambda _{1},\cdots ,\lambda _{q-1}=0}
\frac{(-1)^{\lambda _1+\cdots +\lambda _{q-1}+N}s_1^{\lambda
_1}\cdots s_{q-1}^{\lambda _{q-1}}s_q^Nt^{\lambda _1+\cdots +\lambda
_{q-1}+N}}{\lambda _{1}!\cdots\lambda
_{q-1}!}{\rm Str}[f^0[c(df^1)]^{[\lambda _{1}]}\right.$$
$$\left.\cdots[c(df^{q-1})]^{[\lambda _{q-1}]}\{ [c(df^q)]^{[N]}(s_qt)\}
e^{-(s_{q+1}-s_q)tH_X}\cdots
c(df^{2k})e^{-(1-s_{2k})tH_X}]_Jdv\right|$$
$$\sim O(t^{\frac{2k-2J+\lambda_1+\cdots+\lambda_{q-1}+N-{\rm dim}M}{2}}).\eqno(2.14)$$
Therefore,\\

\noindent{\bf Theorem 2.6}~~ (1)~{\it if~ $2k\leq 2J+{\rm dim}M$, then}
$${{\bf {\rm ch}}_{2k}}(\sqrt{t}D,X)(f^0,\cdots ,
f^{2k})_J~~~~~~~~~~~~~~~~~~~~~~~~~~~~~~~~~~~~~~~~~~~$$
$$=\psi_t
\sum^{{\rm dim}M+2J-2k}_{\lambda_1,\dots,\lambda_{2k}=0}
\frac{
(-1)^{\lambda_1+\cdots+\lambda_{2k}}}{\lambda_1!\cdots \lambda _{2k}!}Ct^{|\lambda|+{k}}
{\rm
Str}[f^0[c(df^1)]^{[\lambda_1]}\cdots [c(df^{2k})]^{[\lambda_{2k}]}e^{-tH_X}]_J+O(\sqrt {t}),\eqno(2.15)$$
{\it with the constant}\\
$$C=\frac{1}{\lambda_1+1}\frac{1}{\lambda_1+\lambda_2+2}\cdots \frac{1}{\lambda_1+\cdots +\lambda_{2k}+2k}.\eqno(2.16)$$
(2) {\it if~ $2k>2J+{\rm dim}M$, then}
$${{\bf {\rm ch}}_{2k}}(\sqrt{t}D,X)(f^0,\cdots ,
f^{2k})_J=O(\sqrt {t}).\eqno(2.17)$$\\

\noindent{\bf 2.2 Computations of infinitesimal equivariant Chern-Connes characters}\\

 Since $H_X$ is a generalized Laplacian, the heat operator $e^{-tH_X}$ exists and
 $$(\frac{\partial}{\partial t}+H_X)e^{-t{H_X}}=0,~~H_Xe^{-t{H_X}}=e^{-t{H_X}}H_X.\eqno(2.18)$$
 It is easy to extend the notation of the Volterra pseudodifferential operator to the case with coefficients in ${\mathbb{C}}[{\mathfrak{g}}^*]$
  (see [BGS],[Gr],[Po]). Let $Q=(H_X+\frac{\partial}{\partial t})^{-1}$ be the Volterra inverse of $H_X+\frac{\partial}{\partial t}$ as in [BGS]. Let
 $K_{Q}(x,y,X,t)$, $k(x,y,X,t)$ be the distribution kernel of $Q$ and the heat kernel of $e^{-tH_X}$ respectively. Then for $t>0$ (see [BGS])
 $$k(x,y,X,t)=K_{Q}(x,y,X,t)+O(t^{\infty})~~~~~~~~~~{\rm as}~~t\rightarrow 0^+.\eqno(2.19)$$
 For the definition 2.4 in [Wa2], we replace $\wedge T_z^* B$ by ${\mathbb{C}}[{\mathfrak{g}}^*]$
  so that we can define Volterra symbols with coefficients in ${\mathbb{C}}[{\mathfrak{g}}^*]$
   and Volterra pseudodifferential operators with coefficients in ${\mathbb{C}}[{\mathfrak{g}}^*]$.
   We denote the space of Volterra pseudodifferential operators with coefficients in ${\mathbb{C}}[{\mathfrak{g}}^*]$ by
  $\Psi_V^*({\mathbb{R}}^n\times {\mathbb{R}}, S(TM)\otimes {\mathbb{C}}[{\mathfrak{g}}^*]).$ \\
  \indent Recall that the quantization map $c:\wedge
 T^*_{{\mathbb{C}}}(M)\rightarrow {\rm Cl}(M)$ and the symbol map
 $\sigma=c^{-1}$ satisfy
 $$\sigma(c(\xi)c(\eta))=\xi\wedge\eta-\left<\xi,\eta\right>.\eqno(2.20)$$
Thus, for $\xi$ and $\eta$ in $\wedge T^*_{{\mathbb{C}}}(M)$ we have
$$\sigma(c(\xi^{(i)})c(\eta^{(j)}))=\xi^{(i)}\wedge \eta^{(j)}~~{\rm
mod}~~\wedge^{i+j-2} T^*_{{\mathbb{C}}}(M),\eqno(2.21)$$ where
$\xi^{(l)}$ denotes the component in $\wedge^{l}
T^*_{{\mathbb{C}}}(M)$ of $\xi\in \wedge T^*_{{\mathbb{C}}}(M).$
Recall that if $e_1,\cdots,e_n$ is an orthonormal frame of $T_xM$,
then
$${\rm Str}[c(e^{i_1})\cdots c(e^{i_k})]=\left\{\begin{array}{lcr}0~& {\rm when}~ k<n,
\\
(-2i)^{\frac{n}{2}} & {\rm when }~ k=n.
\end{array}\right.
\eqno(2.22)$$
 We compute the Chern-Connes character
at a fixed point $x_0\in M$. Using normal coordinates
centered at $x_0$ in $M$ and paralleling $\partial_i$ at $x_0$
along geodesics through $x_0$, we get the orthonormal frame
$e_1,\cdots,e_n$.  We define the Getzler order as follows:
 $${\rm deg}\partial_j=\frac{1}{2}{\rm deg}\partial_t={\rm deg}c(dx_j)=\frac{1}{2}{\rm deg}(X)=-{\rm deg}x^j=1.\eqno(2.23)$$
Let $Q\in \Psi_V^*({\mathbb{R}}^n\times {\mathbb{R}}, S(TM)\otimes
{\mathbb{C}}[{\mathfrak{g}}^*])$ have the symbol $$q(x,X,\xi,\tau)\sim
 \sum_{k\leq m'}q_{k}(x,X,\xi,\tau),\eqno(2.24)$$
 where $q_{k}(x,X,\xi,\tau)$ is an order $k$ symbol. Then using Taylor expansions at $x = 0$ as well as at $X=0$, it gives that
$$\sigma[q(x,X,\xi,\tau)]\sim\sum_{j,k,\alpha,\beta}\frac{x^\alpha}{\alpha!}\frac{X^\beta}{\beta!}
\sigma[\partial_x^\alpha \partial_X^\beta
q_{k}(0,0,\xi,\tau)]^{(j)}.\eqno(2.25)$$ The symbol
$\frac{x^\alpha}{\alpha!}\frac{X^\beta}{\beta!}
\sigma[\partial_x^\alpha \partial_X^\beta
q_{k}(0,0,\xi,\tau)]^{(j)}$ is the Getzler homogeneous
of $k+j-|\alpha|+2|\beta|$.\\

 \noindent {\bf Definition 2.7} The $J$-truncated symbol of $q$ is defined by
 $$\sigma[q(x,X,\xi,\tau)]_J:=\sum_{j,k,\alpha,|\beta|\leq J}\frac{x^\alpha}{\alpha!}\frac{X^\beta}{\beta!}
\sigma[\partial_x^\alpha \partial_X^\beta
q_{k}(0,0,\xi,\tau)]^{(j)}.\eqno(2.26)$$

 Then $\sigma[q(x,X,\xi,\tau)]_J$ can be written as
$$\sigma[q(x,X,\xi,\tau)]_J\sim \sum_{l\geq 0}q_{(m-l)}(x,X,\xi,\tau)_J,~~~~~~~~~q_{(m)}\neq 0, \eqno(2.27)$$
where $q_{(m-l),J}$ is a Getzler homogeneous symbol of degree $m-l$, and the degree of $X$ is $\leq J$.\\

\noindent {\bf Definition 2.8} The integer $m$ is called the Getzler order of $Q$. The symbol $q_{(m),J}$ is the truncated principle Getzler
homogeneous symbol of $Q$. The operator $Q_{(m),J}=q_{(m)}(x,D_x,D_t)_J$ (see [BGS], [Po]) is called the truncated model operator of $Q$.\\

\noindent {\bf Lemma 2.9} {\it Let
$Q\in\Psi_V^*({\mathbb{R}}^n\times {\mathbb{R}}, S(TM)\otimes
{\mathbb{C}}[{\mathfrak{g}}]^*),$ and $Q_J$ has the Getzler order
$m$ and the model operator $Q_{(m),J}$. Then as $t\rightarrow 0^+,$ we
have:}
$${\rm 1)}~~ \sigma[K_{Q_J}(0,0,\frac{X}{t},t)]^{(j)}=O(t^{\frac{j-n-m-1}{2}}),~~{\rm if~} m-j~~{\rm ~is~ odd};~~~~~~~~~~~~~~~~~~~~~~~~~~~~~~~~~~~$$
$${\rm 2)}~~\sigma[K_{Q_J}(0,0,\frac{X}{t},t)]^{(j)}=t^{\frac{j-n-m-2}{2}}K_{Q_{(m),J}}(0,0,X,1)^{(j)}+O(t^{\frac{j-n-m}{2}}),~~{\rm if~} m-j~~{\rm ~is~ even},$$
{\it where $[K_{Q_J}(0,0,\frac{X}{t},t)]^{(j)}$ denotes taking the
$j$ degree form component in $\wedge^*T^*M$. In particular, when $m=-2$ and $j=n$ is
even, we get}
$$\sigma[K_{Q_J}(0,0,\frac{X}{t},t)]^{(n)}=K_{Q_{(-2),J}}(0,0,X,1)^{(n)}+O(t).\eqno(2.28)$$\\

\noindent {\bf Proof.} By (1.7) in [Po], we have
$$K_{Q_J}(0,0,\frac{X}{t},t)\sim\sum_{m_0-j_0~{\rm even}}t^{\frac{j_0-n-m_0-2}{2}}\check{q}_{m_0-j_0}(0,0,\frac{X}{t},1)_J,\eqno(2.29)$$
where $m_0$ is the operator order of $Q_J$. Then
$$\sigma[K_{Q_J}(0,0,\frac{X}{t},t)]^{(j)}\sim\sum_{m_0-j_0~{\rm even}}\sum_{|\beta|\leq J}t^{\frac{j_0-n-m_0-2|\beta|-2}{2}}\sigma[
\frac{X^\beta}{\beta!}\frac{\partial}{\partial
X^\beta}\check{q}_{m_0-j_0}(0,0,0,1)]_J^{(j)}.\eqno(2.30)$$ Let
$L=m_0-j_0+j+2|\beta|$. By $Q_J$ having the Getzler order $m$, then
$L\leq m$.
Thus
$$\sigma[K_{Q_J}(0,0,\frac{X}{t},t)]^{(j)}\sim\sum_{m_0-j_0~{\rm
even}}\sum_{|\beta|\leq J}t^{\frac{j-n-L-2}{2}}\sigma[
\frac{X^\beta}{\beta!}\frac{\partial}{\partial
X^\beta}\check{q}_{m_0-j_0}(0,0,0,1)]_J^{(j)}.\eqno(2.31)$$
 Note that the degree of the leading term is $L=m$ and
$m_0-j_0=m-j-2|\beta|$. When $m-j$ is odd, as $m_0-j_0$ is even,
it is impossible. Therefore,
$$ \sigma[\psi_tK_{Q_J}(0,0,\frac{X}{t},t)]^{(j)}=O(t^{\frac{j-n-m-1}{2}}).\eqno(2.32)$$
 When $L=m$ and $m-j$ is even, the leading coefficient is
$$\sigma[\breve{q}_{(m)}(0,0,X,1)]_J^{(j)}=\sum_{|\beta|\leq J}\sigma[\frac{X^\beta}{\beta!}\frac{\partial}{\partial
X^\beta}\breve{q}_{m-j}(0,0,0,1)]^{(j)}=K_{Q_{(m),J}}(0,0,X,1)^{(j)}.\eqno(2.33)$$
For the next term, it is that $L=m-1$, $m-j$ is even, $m_0-j_0+j=m-1$,
which is impossible, so that the next term is $O(t^{\frac{j-n-m}{2}})$.
~~~~~$\Box$\\

\indent Let $\theta_X$ be the one-form associated with $X_M$ which is defined
by $\theta_X(Y)=g(X,Y)$ for the vector field $Y$. Let $\nabla^{S,X}$
be the Clifford connection $\nabla^S-\frac{1}{4}\theta_X$ on the
spinors bundle and $\triangle_X$ be the Laplacian on $S(TM)$
associated with $\nabla^{S,X}$. Let
$\mu(X)(\cdot)=\nabla^{TM}_{\cdot}X_M$. Define $\alpha:U\times
{\mathfrak{g}}\rightarrow {\mathbb{C}}$ via the formula
$$\alpha_X(x):=-\frac{1}{4}\int_0^1(\iota({\mathcal{R}})\theta_X)(tx)t^{-1}dt,~~\rho(X,x)=e^{\alpha_X(x)},\eqno(2.34)$$
where ${\mathcal{R}}=\sum_{i=1}^nx_i\frac{\partial}{\partial_{x_i}}.$ Recall\\

\noindent{\bf Lemma 2.10} ([BGV Lemma 8.13]) {\it ~The following identity holds}
$$H_X=\Delta_X+\frac{1}{4}r_M,\eqno(2.35)$$
 {\it where $r_M$ is the scalar curvature. In the
trivialization of $S(TM)$ over $U$, the conjugate
$\rho(X,x)(\nabla^{S,X}_{\partial_i})\rho(X,x)^{-1}$ is given by }
$$\rho(X,x)(\nabla^{S,X}_{\partial_i})\rho(X,x)^{-1}=\partial_i-
\frac{1}{4}\sum_{j,a<b}\left<R(\partial_i,\partial_j)e_a,e_b \right>
c(e^a)c(e^b)x^j-\frac{1}{4}\mu_{ij}^M(X)x^j+O_G(0),\eqno(2.36)$$
{\it where $O_G(0)$ is the Getzler order $0$ operator.}\\

By Lemma 2.10, we get \\

\noindent {\bf Proposition 2.11} {\it In the trivialization of
$S(TM)$ over $U$ and the normal coordinate, the model operator of
$\rho(X,x)H_X\rho(X,x)^{-1}$ is}
$$(\rho(X,x)H_X\rho(X,x)^{-1})_{(2)}=-\sum_{i=1}^n(\partial_i-\frac{1}{4}\sum_{j=1}^na_{ij}x_j)^2,~~~~a_{ij}=\left<R^{TM}\partial_i,\partial_j\right>
+\left<\mu(X)\partial_i,\partial_j\right>.\eqno(2.37)$$\\

Let $$\widetilde{[c(df^j)]}^{[\lambda_j]}=[\rho
H_X\rho^{-1},\widetilde{[c(df^j)]}^{[\lambda_j-1]}];~~\widetilde{[c(df^j)]}^{[0]}=c(df^j).$$
Then $$\rho
{[c(df^j)]}^{[\lambda_j]}\rho^{-1}=\widetilde{[c(df^j)]}^{[\lambda_j]};~~O_G(\rho
{[c(df^j)]}^{[\lambda_j]}\rho^{-1})=2\lambda_j,~~{\rm for
}~~\lambda_j>0.\eqno(2.38)$$ We will compute $${\rm
lim}_{t\rightarrow 0}t^{|\lambda|+k}\psi_t{\rm
Str}[f^0[c(df^1)]^{[\lambda_1]}\cdots
[c(df^{2k})]^{[\lambda_{2k}]}e^{-tH_X}]_J.$$ By $\rho e^{-tH_X}
\rho^{-1}=e^{-t \rho H_X \rho^{-1}}$ and (2.38), for a fixed point
$x_0$, then we have
$${\rm
lim}_{t\rightarrow 0}t^{|\lambda|+k}\psi_t{\rm
Str}[f^0[c(df^1)]^{[\lambda_1]}\cdots
[c(df^{2k})]^{[\lambda_{2k}]}e^{-tH_X}]_J$$ $$ ={\rm
lim}_{t\rightarrow 0}t^{|\lambda|+k}\psi_t{\rm
Str}[f^0\widetilde{[c(df^1)]}^{[\lambda_1]}\cdots
\widetilde{[c(df^{2k})]}^{[\lambda_{2k}]}e^{-t\rho
H_X\rho^{-1}}]_J.\eqno(2.39)$$ By (2.38), when
$(\lambda_1,\cdots,\lambda_{2k})\neq (0,\cdots, 0)$, then
$$O_G(f^0\widetilde{[c(df^1)]}^{[\lambda_1]}\cdots
\widetilde{[c(df^{2k})]}^{[\lambda_{2k}]})=O_G(2|\lambda|+2k-1);$$
$$~~O_G(f^0\widetilde{[c(df^1)]}^{[\lambda_1]}\cdots
\widetilde{[c(df^{2k})]}^{[\lambda_{2k}]}(\rho H_X
\rho^{-1}+\partial_t)^{-1})=O_G(2|\lambda|+2k-3).\eqno(2.40)$$ By
(2.40),(2.22) and Lemma 2.9,
$${\rm
lim}_{t\rightarrow 0}t^{|\lambda|+k}\psi_t{\rm
Str}[f^0[c(df^1)]^{[\lambda_1]}\cdots
[c(df^{2k})]^{[\lambda_{2k}]}e^{-tH_X}]_J=0.\eqno(2.41)$$
 When
$(\lambda_1,\cdots,\lambda_{2k})=(0,\cdots, 0),$ then
$O_G(f^0c(df^1)\cdots c(df^{2k}))=2k$
 and
$$O_G(f^0c(df^1)\cdots c(df^{2k})(\rho H_X
\rho^{-1}+\partial_t)^{-1})=O_G(2k-2).$$ The model operator of
$f^0c(df^1)\cdots c(df^{2k})(\rho H_X \rho^{-1}+\partial_t)^{-1}$ is
$$f^0\wedge df^1\wedge\cdots\wedge df^{2k}((\rho H_X
\rho^{-1})_{(2)}+\partial_t)^{-1}.$$ By Lemma 2.9 and Proposition 2.11 in connection with the Mehler formula, we get
$${\rm lim}_{t\rightarrow 0^+}t^k\sigma[\psi_tf^0c(df^1)\cdots c(df^{2k})e^{-tH_X}]_J^{(n)} ={(2\pi\sqrt{-1})}^{-n/2}[f^0\wedge
df^1\wedge\cdots\wedge
df^{2k}\widehat{A}(F^M_{\mathfrak{g}}(X))]^{(n)}_J,\eqno(2.42)$$ where
$\widehat{A}(F^M_{\mathfrak{g}}(X))$ is the equivariant
$\widehat{A}$-genus. By (2.41), (2.42) and Theorem 2.6, we get when $J\rightarrow +\infty$ that\\

\noindent{\bf Theorem 2.12} {\it When $2k\leq {\rm dim}M$ and $X$ is
small which means that $||X_M||$ is sufficient small, then for $f^j\in C^{\infty}_G(M)$,}
$${\rm lim}_{J\rightarrow +\infty}{\rm lim}_{t\rightarrow 0}{{\bf {\rm ch}}_{2k}}(\sqrt{t}D,X)(f^0,\cdots ,
f^{2k})_J$$
$$=\frac{1}{(2k)!}{(2\pi\sqrt{-1})}^{-n/2}\int_Mf^0\wedge
df^1\wedge\cdots\wedge
df^{2k}\widehat{A}(F^M_{\mathfrak{g}}(X))d{\rm
Vol}_M.\eqno(2.43)$$\\

\noindent{\bf Remark.} Theorem 2.12 is not direct from the
equivariant Chern-Connes character formula due to Chern-Hu in [CH]
and the localization formula because $f^0\wedge df^1\wedge\cdots\wedge
df^{2k}\widehat{A}(F^M_{\mathfrak{g}}(X))$ is not an equivariant closed form.\\

Let $p\in M_r({\mathbb{C}}^{\infty}(M))$ be a selfadjoint idempotent, and
$${\rm Ch}({\rm
Im}(p))=\sum_{k=0}^{\infty}(-\frac{1}{2\pi\sqrt{-1}})^k\frac{1}{k!}{\rm
Tr}[p(dp)^{2k}].\eqno (2.44)$$
Let $D_{{\rm Im}p}$ be the Dirac operator with coefficients from ${\rm Im}p.$ Let $S(TM)=S^+(TM)\oplus S^-(TM)$ and $D_{{\rm Im}p,+}$ be the restriction on
$S^+(TM)\otimes {\rm Im}p.$
Then, by the infinitesimal equivariant index formula and Theorem 2.12, we get\\

\noindent{\bf Corollary 2.13} {\it When $X$ is
small, we have}
$${\rm Ind}_{e^{-X}}(D_{{\rm Im}p,+})={\rm lim}_{J\rightarrow +\infty}{\rm lim}_{t\rightarrow 0}\left<{\bf {\rm ch}}_{{\rm even}}
 (\sqrt{t}D,X)_J, {\rm ch}(p)\right>.\eqno(2.45)$$\\

    Next, we shall give an infinitesimal equivariant index formula
for odd dimensional manifolds. Let $M$ be a compact oriented odd dimensional Riemannian manifold
without boundary with a fixed spin structure and $S$ be the bundle
of spinors on $M$. The fundamental setup consistents with that in Section
2.1. Let $g\in GL_r({{C}}^{\infty}(M))$, $g(hx)=g(x)$ for $h\in G$ and $x\in M$.  For $0\leq u\leq 1$, on the bundle $S(TM)\otimes C^r$, let
$$D_{-X,u}=(1-u)D_{-X}+ug^{-1}D_{-X}g=D_{-X}+ug^{-1}dg,~~A=g^{-1}dg,\eqno(2.46)$$
$$H_{X,u}=D_{-X,u}^2+uA(X_M)+L_X.\eqno(2.47)$$
We will compute
$${\rm lim}_{t\rightarrow 0}\int_0^1\sqrt{t}{\rm Tr}\left[Ae^{-tH_{\frac{X}{t},u}}\right]du.$$
By Lemma 2.10, we have\\

\noindent {\bf Proposition 2.14}{\it~ The following identity holds}
$$H_{X,u}=\triangle_X+\frac{1}{4}r_M+u^2c(A)^2+u(D^{cl}(c(A))-2\nabla_A^{S(TM),X}),\eqno(2.48)$$
{\it where $D^{cl}$ is the Dirac operator on the Clifford bundle.}\\

By Lemma 2.10 and Proposition 2.11, we get\\

\noindent {\bf Proposition 2.15} {\it In the trivialization of
$S(TM)$ over $U$ and the normal coordinate, the model operator of
$\rho(X,x)H_{X,u}\rho(X,x)^{-1}$ is}
$$(\rho(X,x)H_{X,u}\rho(X,x)^{-1})_{(2)}=-\sum_{i=1}^n(\partial_i-\frac{1}{4}\sum_{j=1}^na_{ij}x_j)^2+u^2
A^2+udA,\eqno(2.49)$$
$$~~~~a_{ij}=\left<R^{TM}\partial_i,\partial_j\right>
+\left<\mu(X)\partial_i,\partial_j\right>.\eqno(2.50)$$\\

By Lemma 2.9 2), similarly to Theorem 2.12, we get\\

\noindent{\bf Theorem 2.16} {\it When $X$ is
small, then}
$${\rm lim}_{t\rightarrow 0}\int_0^1\sqrt{t}{\rm Tr}\left[Ae^{-tH_{\frac{X}{t},u}}\right]du
={(2\pi\sqrt{-1})}^{-n/2}\int_M\widehat{A}(F^M_{\mathfrak{g}}(X)){\rm ch}(g)d{\rm
Vol}_M,\eqno(2.51)$$
{\it where the odd Chern character is defined by}
$${\rm ch}(g)=\sum_{k=0}^{+\infty}(-1)^k\frac{k!}{(2k+1)!}{\rm Tr}[(g^{-1}dg)^{2k+1}].\eqno(2.52)$$\\

By Lemma 2.9, we know that Theorem 2.12 also holds for odd dimensional manifolds. So by Theorem 2.16,
we get\\

\noindent{\bf Corollary 2.17} {\it When $X$ is
small, we have}
$${\rm lim}_{t\rightarrow 0}\int_0^1\sqrt{t}{\rm Tr}\left[Ae^{-tH_{\frac{X}{t},u}}\right]du
={\rm lim}_{J\rightarrow +\infty}{\rm lim}_{t\rightarrow 0}\left<{\bf {\rm ch}}_{{\rm odd}}
 (\sqrt{t}D,X)_J, {\rm ch}(g)\right>.\eqno(2.53)$$\\

\section{ Infinitesimal equivariant eta cochains}

\quad Let $N$ be a compact oriented odd dimensional Riemannian
manifold without boundary with a fixed spin structure and $S$ be the
bundle of spinors on $N$. The fundamental setup consistents with that in Section
2.1. Define\\
$${{\bf {\rm ch}}_{k}}(\sqrt{t}D_{-X},D_X)(f^0,\cdots ,
f^{k})_J: =t^{k/2}\sum_{j=0}^k\int_{\triangle_{k+1}}{\rm
Tr}\left[\psi_t e^{-t{
\mathfrak{L}}_X}f^0e^{-s_1t(D+\frac{1}{4}c(X))^2}c(df^1)\right.$$
$$\cdot e^{-(s_2-s_1)t(D+\frac{1}{4}c(X))^2}\cdots
c(df^j)e^{-(s_{j+1}-s_j)t(D+\frac{1}{4}c(X))^2}$$
$$\left.\cdot D_Xe^{-(s_{j+2}-s_{j+1})t(D+\frac{1}{4}c(X))^2}c(df^{j+1})\cdots
 c(df^{k})e^{-(1-s_{k+1})t(D+\frac{1}{4}c(X))^2}\right]_Jd{\rm Vol}_{\Delta_{k+1}},\eqno(3.1)$$
 where
$\triangle_{k+1}=\{(s_1,\cdots,s_{k+1})|0\leq s_1\leq
s_2\leq\cdots\leq s_{k+1}\leq 1\}$ is the $k+1$-simplex. Formally,
{\bf truncated infinitesimal equivariant $\eta$ cochains} on
$C_G^{\infty}(N)$ are defined by formulas:\\
$$\widetilde{\eta}_{X,k}(D)_J=\frac{1}{\Gamma(\frac{1}{2})}\int^{\infty}_{\varepsilon}\frac{1}{2\sqrt{t}}{\bf {\rm ch}}_k(\sqrt{t}D_{-X},D_X)_Jdt,
\eqno(3.2)$$
$$\eta_{X,k}(D)_J=\frac{1}{\Gamma(\frac{1}{2})}\int^{\infty}_{\varepsilon}\frac{1}{2\sqrt{t}}{\bf {\rm ch}}_k(\sqrt{t}D_{-X},D_{-X})_Jdt,
\eqno(3.3)$$ \noindent where $\Gamma(\frac{1}{2})=\sqrt {\pi}$ and
$\varepsilon$ is a small positive number. Then
$\widetilde{\eta}_{X,0}(D)(1)$ is the half of the truncated
infinitesimal equivariant eta invariant defined by Goette in [Go].
In order to prove that the above expression is well defined, it is
necessary to check the integrality near infinity of the integration.
In fact, when $k>{\rm
dim}N+1+2J$, we can take $\varepsilon=0$. First, we prove the regularity at zero.\\

\noindent{\bf Lemma 3.1} {\it When $t\rightarrow 0^+$, then for small $X$ and $f^0,\cdots,f^k\in C_G^{\infty}(N),$ we have}
$$ {{\bf {\rm ch}}_{k}}(\sqrt{t}D_{-X},D_X)_J(f^0,\cdots ,
f^{k})=O(t^{-\frac{1}{2}}).\eqno(3.4)$$ {\it When $k>{\rm
dim}N+1+2J$,}
$$ {{\bf {\rm ch}}_{k}}(\sqrt{t}D_{-X},D_X)_J(f^0,\cdots ,
f^{k})=O(t^{\frac{1}{2}}).\eqno(3.5)$$\\

 \indent In (3.1), the difference between infinitesimal equivariant eta
cochains and equivariant eta cochains is that $D_{-X}$ does not
commute with $D_X$. So we can not apply the trick in [Wa1] directly.
By ${\rm Tr}(AB)={\rm Tr}(BA)$, we have
$$
{\rm Tr}\left[\psi_t e^{-t{
\mathfrak{L}}_X}f^0e^{-s_1t(D+\frac{1}{4}c(X))^2}c(df^1)e^{-(s_2-s_1)t(D+\frac{1}{4}c(X))^2}\cdots
c(df^j)e^{-(s_{j+1}-s_j)t(D+\frac{1}{4}c(X))^2}\right.$$
$$\left.\cdot D_Xe^{-(s_{j+2}-s_{j+1})t(D+\frac{1}{4}c(X))^2}c(df^{j+1})\cdots
 c(df^{k})e^{-(1-s_{k+1})t(D+\frac{1}{4}c(X))^2}\right]_J$$
 $$
={\rm Tr}\left[\psi_t D_Xe^{-(s_{j+2}-s_{j+1})tH_X}c(df^{j+1})\cdots
 c(df^{k})e^{-(1-s_{k+1})tH_X}\right.$$
$$\left.\cdot f^0e^{-s_1tH_X}c(df^1)\cdot e^{-(s_2-s_1)tH_X}\cdots
c(df^j)e^{-(s_{j+1}-s_j)tH_X}\right]_J.\eqno(3.6)$$ By Lemma 2.5, we
commute $e^{-(s_{j+1}-s_j)tH_X}$ with $c(df^j)$ and then commute
heat operators from the right to the left. We write the result for
the case that $k=2,~j=1$. For general case the result is similar.
$$D_Xe^{-(s_{3}-s_{2})tH_X}c(df^{j+1})
 c(df^{2})e^{-(1-s_{3})tH_X}f^0e^{-s_1tH_X}c(df^1)\cdot
 e^{-(s_2-s_1)tH_X}$$
 $$=\sum_{\lambda_1,\lambda_2,\lambda_3=0}^{N-1}\frac{t^{\lambda_1+\lambda_2+\lambda_3}}{\lambda_1!\lambda_2!\lambda_3!}
 (s_2-s_1)^{\lambda_1}s_2^{\lambda_2}
 (1-s_3+s_2)^{\lambda_3}D_Xe^{-tH_X}c(df^2)^{[\lambda_3]}(f^0)^{[\lambda_2]}c(df^1)^{[\lambda_1]}$$
$$
+\sum_{\lambda_1,\lambda_2=0}^{N-1}\frac{t^{\lambda_1+\lambda_2+N}}{\lambda_1!\lambda_2!}(s_2-s_1)^{\lambda_1}s_2^{\lambda_2}
 (1-s_3+s_2)^{N}D_X$$
 $$\cdot e^{-t(s_3-s_2)H_X}\{c(df^2)_1^{[N]}[(1-s_3+s_2)t]\}(f^0)^{[\lambda_2]}c(df^1)^{[\lambda_1]}$$
$$+\sum_{\lambda_1=0}^{N-1}\frac{t^{\lambda_1+N}}{\lambda_1!}(s_2-s_1)^{\lambda_1}s_2^{N}
 D_Xe^{-t(s_3-s_2)H_X}c(df^2)e^{-t(1-s_3)H_X}[(f^0)_1^{[N]}(t_1s_2)]c(df^1)^{[\lambda_1]}$$
$$+t^{N}(s_2-s_1)^{N}
 D_Xe^{-t(s_3-s_2)H_X}c(df^2)e^{-t(1-s_3)H_X}f^0e^{-ts_1H_X}[c(df^1)_1^{[N]}[(s_2-s_1)t]].\eqno(3.7)$$
 For the second term on the right hand side of (3.7), we have
 $$A:=t^{k/2}\int_{\triangle_{k+1}}\left|{\rm Tr}\left[\psi_t
 \sum_{\lambda_1,\lambda_2=0}^{N-1}\frac{t^{\lambda_1+\lambda_2+N}}{\lambda_1!\lambda_2!}(s_2-s_1)^{\lambda_1}s_2^{\lambda_2}
 (1-s_3+s_2)^{N}D_Xe^{-t(s_3-s_2)H_X}\right.\right.$$
 $$\left.\left.\cdot\int_{\triangle_N}e^{-t(1-u_1)(1-s_3+s_2)H_X}c(df^2)^{[N]}e^{-tu_1(1-s_3+s_2)H_X}(f^0)^{[\lambda_2]}c(df^1)
 ^{[\lambda_1]}\right]_J\right|du_1\cdots du_nds_1ds_2ds_3
 $$
 $$=\sum_{\lambda_1,\lambda_2=0}^{N-1}\frac{t^{\lambda_1+\lambda_2+N+k/2-J}}{\lambda_1!\lambda_2!}\int_{\triangle_{k+1}}\int_{\triangle_N}
 (s_2-s_1)^{\lambda_1}s_2^{\lambda_2}
 (1-s_3+s_2)^{N}$$
 $$\cdot \left|{\rm Tr}\left[D_Xe^{-t\sigma_1H_X}e^{-t\sigma_2H_X}c(df^2)^{[N]}e^{-t\sigma_3H_X}(f^0)^{[\lambda_2]}c(df^1)
 ^{[\lambda_1]}\right]_J\right|du_1\cdots du_nds_1ds_2ds_3,\eqno(3.8)$$
 where $\sigma_1+\sigma_2+\sigma_3=1,~\sigma_1,\sigma_2,\sigma_3\geq 0$ and
 $$\sigma_1=s_3-s_2;~~\sigma_2=(1-u_1)(1-s_3+s_2),~~\sigma_3=u_1(1-s_3+s_2).\eqno(3.9)$$
 We divide the region into three parts as shown in Lemma 2.2. By the Weyl theorem, we get that when $N\geq n+2-k+2J$, then
 $$A\sim O(t^{\frac{N+k+|\lambda|-n-1-2J}{2}})\sim O(t^{\frac{1}{2}}).\eqno(3.10)$$
 Similarly, we get that when $N\geq n+2-k+2J$, the third and fourth terms on the right hand side in (3.7) are also $O(t^{\frac{1}{2}})$. When $k\geq n+2+2J$, then
 $N\geq n+2-k+2J$. So we get\\

\noindent{\bf Theorem 3.2}~~{\it {\rm 1)}~~If $k\leq n+1+2J$ and $X$ is small,
then
when $t\rightarrow 0^+$, we have:}\\
$${\bf {\rm ch}}_{k}(\sqrt
{t}D_{-X},D_X)(f^0,\cdots ,f^k)_J
=\sum^k_{j=0}(-1)^j
\sum_{0\leq \lambda _{1},\cdots ,\lambda _k\leq
{N-1}}\frac{(-1)^{|\lambda|}C't^{|\lambda|+\frac{k}{2}}}{\lambda!}$$
$$\cdot
{\rm Tr}\left[\psi_tc(df^{j+1})^{[\lambda_{k+1}]}\cdots
c(df^{k})^{[\lambda_{j+2}]}(f^0)^{[\lambda_{j+1}]}c(df^{1})^{[\lambda_{j}]}\cdots
c(df^{j})^{[\lambda_{1}]}D_Xe^{-tH_X}\right]_J
+O(t^{\frac{1}{2}}),\eqno(3.11)$$
{\it where $C'$ is a constant.}\\
{\it {\rm 2)}~~If $k>n+1+2J$, then
when $t\rightarrow 0^+$, we have:}\\
$${\bf {\rm ch}}_{k}(\sqrt
{t}D_{-X},D_X)(f^0,\cdots ,f^k)_J\sim O(t^{\frac{1}{2}}).\eqno(3.12)$$\\
{\it 1),2) also hold for ${\bf {\rm ch}}_{k}(\sqrt
{t}D_{-X},D_{-X})(f^0,\cdots ,f^k)_J.$ }\\

 \noindent{\bf Lemma 3.3}~~{\it When
$t\rightarrow 0^+$, we have:}
$$t^{|\lambda|+\frac{k}{2}}{\rm
Tr}\left[\psi_tc(df^{j+1})^{[\lambda_{k+1}]}\cdots
c(df^{k})^{[\lambda_{j+2}]}\right.$$
$$\left.\cdot (f^0)^{[\lambda_{j+1}]}c(df^{1})^{[\lambda_{j}]}\cdots c(df^{j})^{[\lambda_{1}]}D_Xe^{-tH_X}\right]_J\sim O(t^{-\frac{1}{2}});\eqno(3.13)$$
$$t^{|\lambda|+\frac{k}{2}}\int_0^1{\rm
Tr}\left[\psi_tc(df^{j+1})^{[\lambda_{k+1}]}\cdots
c(df^{k})^{[\lambda_{j+2}]}\right.$$
$$\left.\cdot (f^0)^{[\lambda_{j+1}]}c(df^{1})^{[\lambda_{j}]}\cdots c(df^{j})^{[\lambda_{1}]}e^{-t\sigma_0H_X}D_X
e^{-t(1-\sigma_0)H_X}\right]_Jd\sigma_0\sim
O(t^{\frac{1}{2}}).\eqno(3.14)$$\\

\noindent{\bf Proof.} We introduce an
auxiliary Grassmann variable $z$ as shown in [BF]. Let
$$\widetilde{H_X}=H_X-zD_X;~h(x)=1+\frac{1}{2}z\sum_{j=1}^nx_ic(e_i).\eqno(3.15)$$
Then we have by Lemma 2.10 that
$$\widetilde{H_X}=-\sum_{j=1}^n(\nabla^{S,X}_{e_j}-\frac{1}{2}c(e_j)z)^2+\sum_{j=1}^n(\nabla^{S,X}_{\nabla^{TM}_{e_j}e_j}
-\frac{1}{2}c(\nabla^{TM}_{e_j}e_j)z)+\frac{1}{4}r_M.\eqno(3.16)$$
Using Lemma 8.13 in [BGV], we have
$$\rho\widetilde{H_X}\rho^{-1}=-\sum_{j=1}^n(\nabla^{S}_{e_j}-\frac{1}{4}\sum_i\left<\mu_X(e_j),\partial_i\right>x^i
+\left<{h_j}(x),X\right>
-\frac{1}{2}c(e_j)z)^2$$
$$+\sum_{j=1}^n(\nabla^{S}_{\nabla^{TM}_{e_j}e_j}
-\frac{1}{2}c(\nabla^{TM}_{e_j}e_j)z-\frac{1}{4}\sum_i\left<\mu_X(\nabla^{TM}_{e_j}e_j),\partial_i\right>x^i
+\left<\overline{h_j}(x),X\right>)+\frac{1}{4}r_M,\eqno(3.17)$$
where ${h_j}(x), ~\overline{h_j}(x)=O(|x|^2).$ Then
$$(h\rho)\widetilde{H_X}(h\rho)^{-1}=\rho H_X\rho^{-1}+zu,~~{\rm where}~ O_G(u)\leq 0 {\rm ~has ~no } ~z.\eqno(3.18)$$
By the Duhamel principle, we have
$${\rm exp}(-t\widetilde{H}_X)={\rm exp}(-t{H}_X)+tz\int_0^1e^{-t\sigma_0H_X}D_X
e^{-t(1-\sigma_0)H_X}d\sigma_0.\eqno(3.19)$$ By (3.18) and (3.19),
then
$$(h\rho)^{-1}{\rm exp}(-t(\rho H_X\rho^{-1}+zu))(h\rho)=\rho^{-1}{\rm exp}(-t\rho
{H}_X\rho^{-1})\rho$$ $$+tz\int_0^1e^{-t\sigma_0H_X}D_X
e^{-t(1-\sigma_0)H_X}d\sigma_0.\eqno(3.20)$$ Let
$$A_0:=\widetilde{c(df^{j+1})}^{[\lambda_{k+1}]}\cdots
\widetilde{c(df^{k})}^{[\lambda_{j+2}]}
\widetilde{(f^0)}^{[\lambda_{j+1}]}\widetilde{c(df^{1})}^{[\lambda_{j}]}\cdots
\widetilde{c(df^{j})}^{[\lambda_{1}]},$$
$$A_1=c(df^{j+1})^{[\lambda_{k+1}]}\cdots
c(df^{k})^{[\lambda_{j+2}]}
(f^0)^{[\lambda_{j+1}]}c(df^{1})^{[\lambda_{j}]}\cdots
c(df^{j})^{[\lambda_{1}]}.$$
 $${\rm
Tr}[\psi_tA_0 h^{-1}{\rm exp}(-t(\rho H_X\rho^{-1}+zu))h]={\rm
Tr}[\psi_tA_0{\rm exp}(-t\rho {H}_X\rho^{-1})]$$ $$+tz\int_0^1{\rm
Tr}[\psi_tA_1e^{-t\sigma_0H_X}D_X
e^{-t(1-\sigma_0)H_X}]d\sigma_0.\eqno(3.21)$$
$${\rm Tr}[\psi_tA_0
h^{-1}{\rm exp}(-t(\rho H_X\rho^{-1}+zu))h]={\rm Tr}[\psi_th^{-1}A_0
{\rm exp}(-t(\rho H_X\rho^{-1}+zu))h]$$ $$+{\rm Tr}[\psi_t[A_0,
h^{-1}]{\rm exp}(-t(\rho H_X\rho^{-1}+zu))h].\eqno(3.22)$$ Now
$$t^{|\lambda|+\frac{k}{2}}{\rm Tr}[\psi_t[A_0, h^{-1}]{\rm exp}(-t(\rho
H_X\rho^{-1}+zu))h]_J=O(t^{3/2}).\eqno(3.23)$$ In fact, by direct
computations, then when $\lambda\neq (0,\cdots,0)$, we have
$O_G([A_0,h])=2|\lambda|+k-2$ up to terms $x_jL_jz$, where $L_j$ is
an operator. When $\lambda= (0,\cdots,0)$, $[A_0,h]=\sum_jx_jL_j$
and we fix a point $x_0$, so in this case (3.23) is zero. By Lemma 2.9 1), (3.23) is got. By
$$(\partial_t+\rho
H_X\rho^{-1}+zu)^{-1}=(\partial_t+\rho
H_X\rho^{-1})^{-1}-z(\partial_t+\rho
H_X\rho^{-1})^{-1}u(\partial_t+\rho H_X\rho^{-1})^{-1},\eqno(3.24)$$
we have $$t^{|\lambda|+\frac{k}{2}}{\rm Tr}[\psi_tA_0 {\rm
exp}(-t(\rho H_X\rho^{-1}+zu))]_J-t^{|\lambda|+\frac{k}{2}}{\rm
Tr}[\psi_tA_0{\rm exp}(-t\rho
{H}_X\rho^{-1})]_J=O(t^{3/2}).\eqno(3.25)$$ By (3.21)-(3.25), we get
(3.14). Considering $D_{-X}e^{-t\sigma_0H_X}=e^{-t\sigma_0H_X}D_{-X}$, we get
$$e^{-t\sigma_0H_X}D_X
e^{-t(1-\sigma_0)H_X}=D_Xe^{-tH_X}+\frac{1}{2}c(X)e^{-tH_X}-\frac{1}{2}e^{-t\sigma_0H_X}c(X)
e^{-t(1-\sigma_0)H_X}.\eqno(3.26)$$ Using Lemma 2.4, similarly to
Theorem 2.6, we get
$$t^{|\lambda|+\frac{k}{2}}\int_0^1{\rm
Tr}[\psi_tA_1e^{-t\sigma_0H_X}D_X e^{-t(1-\sigma_0)H_X}]d\sigma_0$$
$$=t^{|\lambda|+\frac{k}{2}}{\rm
Tr}[\psi_tA_1D_X e^{-tH_X}] +\sum_{1\leq l\leq
K_0}t^{|\lambda|+\frac{k}{2}}{\rm
Tr}[\psi_tA_1t^lc(X)^{[l]}e^{-tH_X}]+O(t^{1/2}).\eqno(3.27)$$ Considering
$O_G(X)=2$ and $n$ is odd, we get $$\sum_{1\leq l\leq
K_0}t^{|\lambda|+\frac{k}{2}}{\rm
Tr}[\psi_tA_1t^lc(X)^{[l]}e^{-tH_X}]=O(t^{-1/2}).\eqno(3.28)$$ By
(3.14),(3.27) and (3.28), we get (3.13).~~~~~~~ $\Box$\\

\noindent{\bf Remark.} Lemma 2.12 in [Wa1] is not correct. But
using the trick in (3.23), we can prove the regularity of
equivariant eta chains in [Wa1].\\

\indent Next, we prove the regularity at infinity. Let
${\mathcal{M}}$ be the algebra generated by pseudodifferential
operators and smoothing operators. Let ${\mathcal{N}}$ be the ideal of
all smooth operators in ${\mathcal{M}}$. The algebra
${\mathbb{C}}[{\mathfrak{g^*}}]_J$ possesses a natural filtration
$${\mathbb{C}}[{\mathfrak{g^*}}]_{J,j}:=\frac{({\mathfrak{g^*}})^j{\mathbb{C}}[{\mathfrak{g^*}}]}
{({\mathfrak{g^*}})^{J+1}{\mathbb{C}}[{\mathfrak{g^*}}]}.$$ Let
${\mathcal{M}}_j$ be the algebra generated by differential operators and smoothing operators acting on $\Gamma(S(TN))$ with coefficients in
${\mathbb{C}}[{\mathfrak{g^*}}]_{J,j}$. Let ${\mathcal{N}}_j$ denote the algebra generated by smoothing operators acting on $\Gamma(S(TN))$ with coefficients in
${\mathbb{C}}[{\mathfrak{g^*}}]_{J,j}.$
 The elements of
${\mathbb{C}}[{\mathfrak{g^*}}]_{J,j}$ are nilpotent of order $\leq
J+1$ in ${\mathbb{C}}[{\mathfrak{g^*}}]_{J}$ for $j\geq 1$, so the
elements of ${\mathcal{M}}_j$ and ${\mathcal{N}}_j$ are also
nilpotent of the same order. Note that the subspace
$1+{\mathcal{N}}_j$ of ${\mathcal{M}}$ forms a group with inverse
$(1+K_X)^{-1}=\sum_{j=0}^J(-K_X)^j.$ Let $P_0\in {\mathcal{N}}$ be
the projection onto ${\rm ker}(D)$ and set $P_1:=1-P_0\in
{\mathcal{M}}.$ For any $A_X\in {\rm
End}(\Gamma(S(TN)))\otimes{\mathbb{C}}[{\mathfrak{g^*}}]_{J}$ we
shall write
$$A_X=\left|\begin{array}{lcr}
  \  P_0A_XP_0 &  P_0A_XP_1 \\
    \ P_1A_XP_0  & P_1A_XP_1
\end{array}\right|
\in \left|\begin{array}{lcr}
  \  {\mathcal{N}} &  {\mathcal{N}} \\
    \ {\mathcal{N}}  & {\mathcal{M}}
\end{array}\right|.$$\\

\noindent{\bf Lemma 3.4}(Lemma 2.33 [Go]) {\it There exists
$\gamma_X\in 1+{\mathcal{N}}_1$ which commutes with
${\mathfrak{L}}_X$, such that}
$$\gamma_XD_{-X}^2\gamma_{X}^{-1}=\left|\begin{array}{lcr}
  \  U_X &  0 \\
    \ 0  & V_X
\end{array}\right|\in \left|\begin{array}{lcr}
  \  {\mathcal{N}}_2 &  0 \\
    \ 0  & T+{\mathcal{N}}_2
\end{array}\right|,$$
{\it where $T-D^2\in {\mathcal{M}}_1$ and  $U_X$ has the form
$P_0U_X'P_0$.}\\

By Lemma 3.4, we have
$$t\gamma_{\frac{X}{t}}H_{\frac{X}{t}}\gamma_{\frac{X}{t}}^{-1}=\left|\begin{array}{lcr}
  \ tU_{\frac{X}{t}}+P_0{\mathfrak{L}}_XP_0 &  0 \\
    \ 0  & tV_{\frac{X}{t}}+P_1{\mathfrak{L}}_XP_1
\end{array}\right|,\eqno(3.29)$$

$$\gamma_{\frac{X}{t}}D_{\frac{X}{t}}\gamma_{\frac{X}{t}}^{-1}
=\left|\begin{array}{lcr}
  \ 0 &  0 \\
    \ 0  & D
\end{array}\right|+O(t^{-1}).\eqno(3.30)$$

$$e^{-tH_{\frac{X}{t}}}=\gamma_{\frac{X}{t}}^{-1}
\left|\begin{array}{lcr}
  \ P_0e^{-tU_{\frac{X}{t}}-{\mathfrak{L}}_X}P_0 &  0 \\
    \ 0  & P_1e^{-tV_{\frac{X}{t}}-{\mathfrak{L}}_X}P_1
\end{array}\right|\gamma_{\frac{X}{t}},\eqno(3.31)$$
$$e^{-t\sigma_lH_{\frac{X}{t}}}D_{\frac{X}{t}}e^{-t\sigma_{l+1}H_{\frac{X}{t}}}
=\gamma_{\frac{X}{t}}^{-1} \left|\begin{array}{lcr}
  \ 0 &  0 \\
    \ 0  & P_1e^{-t\sigma_l[V_{\frac{X}{t}}+{\mathfrak{L}}_{\frac{X}{t}}]}
    De^{-t\sigma_{l+1}[V_{\frac{X}{t}}+{\mathfrak{L}}_{\frac{X}{t}}]}P_1
\end{array}\right|\gamma_{\frac{X}{t}}$$
$$+\gamma_{\frac{X}{t}}^{-1} \left|\begin{array}{lcr}
  \ P_0e^{-t\sigma_l[U_{\frac{X}{t}}+{\mathfrak{L}}_{\frac{X}{t}}]}P_0 &  0 \\
    \ 0  & P_1e^{-t\sigma_l[V_{\frac{X}{t}}+{\mathfrak{L}}_{\frac{X}{t}}]}P_1
\end{array}\right|L$$
$$\cdot
\left|\begin{array}{lcr}
  \ P_0e^{-t\sigma_{l+1}[U_{\frac{X}{t}}+{\mathfrak{L}}_{\frac{X}{t}}]}P_0 &  0 \\
    \ 0  & P_1e^{-t\sigma_{l+1}[V_{\frac{X}{t}}+{\mathfrak{L}}_{\frac{X}{t}}]}P_1
\end{array}\right|\gamma_{\frac{X}{t}}O(t^{-1})+O(t^{-1}),\eqno(3.32)$$
where $L$ is a zero order operator. We note that
$\gamma_{\frac{X}{t}}=1+O(t^{-1})S_0$ where $S_0$ is a smoothing
operator and we assume that
$\gamma_{\frac{X}{t}}=1$ temporarily.\\

\noindent{\bf Lemma 3.5} {\it When $t\rightarrow +\infty$, we have:}\\
$${\bf {\rm ch}}_{k}(\sqrt
{t}D_{-X},D_X)(f^0,\cdots ,f^k)_J\sim O(t^{-1}).\eqno(3.33)$$\\
{\it It also holds for ${\bf {\rm ch}}_{k}(\sqrt
{t}D_{-X},D_{-X})(f^0,\cdots ,f^k)_J.$ }\\

\noindent{\bf Proof.} Recall that ${\bf {\rm ch}}_{k}(\sqrt
{t}D_{-X},D_X)(f^0,\cdots ,f^k)_J=\sum_{0\leq j\leq k}(-1)^jT_j,$
where
$$T_j =t^{k/2}\int_{\triangle_{k+1}}{\rm
Tr}\left[
f^0e^{-\sigma_0tH_{\frac{X}{t}}}c(df^1)e^{-\sigma_1tH_{\frac{X}{t}}}
\cdots c(df^j)e^{-\sigma_jtH_{\frac{X}{t}}} D_{\frac{X}{t}}\right.$$
$$\left.e^{-\sigma_{j+1}tH_{\frac{X}{t}}} c(df^{j+1})\cdots
 c(df^{k})e^{-\sigma_{k+1}tH_{\frac{X}{t}}}\right]_J
 d{\rm Vol}\Delta_{k+1}.
\eqno(3.34)$$
By (3.31) and (3.32), we
have
$$T_j=t^{k/2}\int_{\triangle_{k+1}}{\rm
Tr}\left[ f^0 [
P_0e^{-\sigma_0(tU_{\frac{X}{t}}+{\mathfrak{L}}_X)}P_0+
   P_1e^{-\sigma_0(tV_{\frac{X}{t}}+{\mathfrak{L}}_X)}P_1]
c(df^1)\right.$$
$$\cdot[P_0e^{-\sigma_1(tU_{\frac{X}{t}}+{\mathfrak{L}}_X)}P_0+
   P_1e^{-\sigma_1(tV_{\frac{X}{t}}+{\mathfrak{L}}_X)}P_1]
\cdots
c(df^j)\left[P_1e^{-t\sigma_j[V_{\frac{X}{t}}+{\mathfrak{L}}_{\frac{X}{t}}]}
    De^{-t\sigma_{j+1}[V_{\frac{X}{t}}+{\mathfrak{L}}_{\frac{X}{t}}]}P_1\right.$$
$$+ (
P_0e^{-t\sigma_j[U_{\frac{X}{t}}+{\mathfrak{L}}_{\frac{X}{t}}]}P_0
   + P_1e^{-t\sigma_j[V_{\frac{X}{t}}+{\mathfrak{L}}_{\frac{X}{t}}]}P_1)$$
   $$\left.\cdot L(
P_0e^{-t\sigma_{j+1}[U_{\frac{X}{t}}+{\mathfrak{L}}_{\frac{X}{t}}]}P_0
+P_1e^{-t\sigma_{j+1}[V_{\frac{X}{t}}+{\mathfrak{L}}_{\frac{X}{t}}]}P_1)O(t^{-1})\right]$$
    $$\left.\cdot
 c(df^{j+1})\cdots
 c(df^{k})[P_0e^{-\sigma_{k+1}(tU_{\frac{X}{t}}+{\mathfrak{L}}_X)}P_0+
   P_1e^{-\sigma_{k+1}(tV_{\frac{X}{t}}+{\mathfrak{L}}_X)}P_1]
\right]_J
 d{\rm Vol}\Delta_{k+1},
\eqno(3.35)$$ where
$P_0e^{-t\sigma_{l+1}[U_{\frac{X}{t}}+{\mathfrak{L}}_{\frac{X}{t}}]}P_0$
and
$P_1e^{-t\sigma_{l+1}[V_{\frac{X}{t}}+{\mathfrak{L}}_{\frac{X}{t}}]}P_1$
stand for respectively
$$\left|\begin{array}{lcr}
  \ P_0e^{-t\sigma_{l+1}[U_{\frac{X}{t}}+{\mathfrak{L}}_{\frac{X}{t}}]}P_0 &  0 \\
    \ 0  & 0
\end{array}\right|~~{\rm and}~~
\left|\begin{array}{lcr}
  \ 0 &  0 \\
    \ 0  & P_1e^{-t\sigma_{l+1}[V_{\frac{X}{t}}+{\mathfrak{L}}_{\frac{X}{t}}]}P_1
\end{array}\right|.$$
 Let
$$T_j'=t^{k/2}\int_{\triangle_{k+1}}{\rm
Tr}\left[ f^0 [
P_0e^{-\sigma_0(tU_{\frac{X}{t}}+{\mathfrak{L}}_X)}P_0+
   P_1e^{-\sigma_0(tV_{\frac{X}{t}}+{\mathfrak{L}}_X)}P_1]
c(df^1)\right.$$
$$\cdot[P_0e^{-\sigma_1(tU_{\frac{X}{t}}+{\mathfrak{L}}_X)}P_0+
   P_1e^{-\sigma_1(tV_{\frac{X}{t}}+{\mathfrak{L}}_X)}P_1]
\cdots
c(df^j)\left[P_1e^{-t\sigma_j[V_{\frac{X}{t}}+{\mathfrak{L}}_{\frac{X}{t}}]}
    De^{-t\sigma_{j+1}[V_{\frac{X}{t}}+{\mathfrak{L}}_{\frac{X}{t}}]}P_1\right]$$
    $$\left.\cdot
 c(df^{j+1})\cdots
 c(df^{k})[P_0e^{-\sigma_{k+1}(tU_{\frac{X}{t}}+{\mathfrak{L}}_X)}P_0+
   P_1e^{-\sigma_{k+1}(tV_{\frac{X}{t}}+{\mathfrak{L}}_X)}P_1]
\right]_J
 d{\rm Vol}\Delta_{k+1},
\eqno(3.36)$$
$$T_j''=t^{k/2-1}\int_{\triangle_{k+1}}{\rm
Tr}\left[ f^0 [
P_0e^{-\sigma_0(tU_{\frac{X}{t}}+{\mathfrak{L}}_X)}P_0+
   P_1e^{-\sigma_0(tV_{\frac{X}{t}}+{\mathfrak{L}}_X)}P_1]
c(df^1)\right.$$
$$\cdot[P_0e^{-\sigma_1(tU_{\frac{X}{t}}+{\mathfrak{L}}_X)}P_0+
   P_1e^{-\sigma_1(tV_{\frac{X}{t}}+{\mathfrak{L}}_X)}P_1]
\cdots c(df^j)$$
$$\left[(
P_0e^{-t\sigma_j[U_{\frac{X}{t}}+{\mathfrak{L}}_{\frac{X}{t}}]}P_0
   + P_1e^{-t\sigma_j[V_{\frac{X}{t}}+{\mathfrak{L}}_{\frac{X}{t}}]}P_1)L(
P_0e^{-t\sigma_{j+1}[U_{\frac{X}{t}}+{\mathfrak{L}}_{\frac{X}{t}}]}P_0
    +P_1e^{-t\sigma_{j+1}[V_{\frac{X}{t}}+{\mathfrak{L}}_{\frac{X}{t}}]}P_1)\right]$$
    $$\left.\cdot
 c(df^{j+1})\cdots
 c(df^{k})[P_0e^{-\sigma_{k+1}(tU_{\frac{X}{t}}+{\mathfrak{L}}_X)}P_0+
   P_1e^{-\sigma_{k+1}(tV_{\frac{X}{t}}+{\mathfrak{L}}_X)}P_1]
\right]_J
 d{\rm Vol}\Delta_{k+1}.
\eqno(3.37)$$ We estimate (3.36) first. Since $P_0c(df^j)P_0=0$,
only the terms containing no more than $\frac{k}{2}+1$ copies of
$P_0e^{-\sigma_l(tU_{\frac{X}{t}}+{\mathfrak{L}}_X)}P_0$ give a
non-zero contribution. In fact, the term containing no copy of
$P_0e^{-\sigma_l(tU_{\frac{X}{t}}+{\mathfrak{L}}_X)}P_0$ has
exponential decay. Note that
$$V_{\frac{X}{t}}+{\mathfrak{L}}_{\frac{X}{t}}=P_1D^2P_1+F_{\frac{X}{t}},\eqno(3.38)$$
where $F_{\frac{X}{t}}\in {\mathcal{M}}_1$. Similarly to Lemma
2.2, by (3.38) we have that when $X$ is small and $t\rightarrow \infty$,
$$||P_1e^{-ut[V_{\frac{X}{t}}+{\mathfrak{L}}_{\frac{X}{t}}]}_JP_1{B}||_{u^{-1}}\leq C(X)_Ju^{-\frac{l}{2}}t^{-\frac{l}{2}}({\rm
tr}[P_1e^{-\frac{tD^2}{2}}P_1])^u,\eqno(3.39)$$ where $B$ is a
$l$-order operator. Using $||{\rm Tr}(P_1e^{-sD^2}P_1)||\leq
C_0e^{-s\lambda^2},~~{\rm for }~s\geq 1$ and (3.39), similarly to Lemma 1.1 in [Wu], we get the exponential decay. Thus, it is left to
deal with the terms containing at least $\frac{k}{2}+1$ copies of
$P_1e^{-\sigma_l(tV_{\frac{X}{t}}+{\mathfrak{L}}_X)}P_1$, as well as
at least one copy of
$P_0e^{-\sigma_l(tU_{\frac{X}{t}}+{\mathfrak{L}}_X)}P_0$, we may use
the trick in Lemma 2 in [CM2] to prove $T_j'=O(t^{-1}).$
Similarly, we can prove $T_j''=O(t^{-2})$.\\
\indent For the general $\gamma_{\frac{X}{t}}=1+O(t^{-1})S_0$, since $P_0c(df^l)S_0P_0$ and $P_0S_0c(df^l)P_0$ do not equal zero, so the number of copies of
$P_1e^{-\sigma_l(tV_{\frac{X}{t}}+{\mathfrak{L}}_X)}P_1$ in (3.36) may be less than $\frac{k}{2}+1$. But the coefficient of $S_0$ is $O(t^{-1})$. Through the careful observation, we still get that (3.36) is $O(t^{-1})$ and then get Lemma 3.5. ~~$\Box$\\

Now we shall give the convergence of the total truncated
infinitesimal equivariant eta cochain. Let $C_G^1(N)$ be Banach
algebra of once differentiable
 function on $N$ with the norm defined by
 $$||f||_1:={\rm sup}_{x\in N}|f(x)|+{\rm sup}_{x\in N}||df(x)||.$$
 Let $$\phi_{X,J}=\{\phi_{X,J,0},\cdots,\phi_{X,J,2q},\cdots\}$$
be a truncated infinitesimal  equivariant even cochains sequence in
the bar complex of $C^1(N).$ Define
$$||\phi_{X,J,2q}||={\rm sup}_{||f_i||_1\leq 1;~0\leq i\leq
 2q}\{|\phi_{X,J,2q}(f_0,\cdots,f_{2q})|\}.$$

 \noindent {\bf Definition 3.6 }~ The radius of convergence of
 $\phi_{X,J}$ is defined to be that of the power series $\sum
 q!||\phi_{X,J,2q}||z^q.$ The space of cochains sequence with radius
 of convergence $r$ at least larger than zero is denoted by $C^{{\rm
 even},X,J}_r(C_G^1(N))$ (define $C^{{\rm
 odd},X,J}_r(C_G^1(N))$ similarly).\\

 \indent In general, the sequence
$$\eta_{X}(D)_J=\{\cdots,\eta_{X,2q}(D)_J,\eta_{X,2q+2}(D)_J,\cdots\}$$
which is called a total truncated infinitesimal equivariant eta cochain is not an entire cochain. Similarly to Proposition 2.16 in [Wa1], we have\\

\noindent {\bf Proposition 3.7} {\it Suppose that $D$ is invertible
with $\lambda$ the smallest positive eigenvalue of $|D|$ and $X$ is small. Then the
truncated infinitesimal equivariant total eta cochain
$\eta_{X}(D)_J$ has radius of convergence $r$ satisfying the
inequality: $r\geq 4\lambda^2>0$ i.e. $\eta_{X}(D)_J\in C^{{\rm
 even},X,J}_{4\lambda^2}(C_G^1(N)).$}\\

\indent For the idempotent $p\in {\cal M}_r(C^{\infty}(N))$, let
$||dp||=||[D,p]||=\sum_{i,j}||dp_{i,j}||$
where $p_{i,j}~(1\leq i,j\leq r)$ is the entry of $p$. Similarly to Proposition 2.17 in [Wa1], we have \\

\noindent {\bf Proposition 3.8} {\it Assume that $D$ is invertible
with $\lambda$ the smallest positive eigenvalue of $|D|$ and $||dp||<\lambda$
and $X$ is small, then
the pairing $\langle\eta_X(D)_J,Ch(p)\rangle$ is well-defined.}\\

 \indent Next we establish the main theorem in this section. Suppose $D$ is
invertible with $\lambda$ the smallest eigenvalue of $|D|$, and
$p=p^*=p^2\in {\cal M }_r(C^{\infty}_G(N))$ is an idempotent which
satisfies $||dp||<\lambda$. Let
$$p(D\otimes I_r)p:~p(H\otimes {\bf C^r})=L^2(N,S\otimes
p({\bf C^r}))\rightarrow L^2(N,S\otimes p({\bf C^r}))$$ be the Dirac
operator with coefficients from $F=p({\bf C^r})$. Since $p\in {\cal M }_r(C^{\infty}_G(N))$,
 we have $$e^{-X}[p(D\otimes I_r)p]=[p(D\otimes I_r)p]e^{-X}.$$
\indent Let
 $${\bf D_{-X}}=\left[\begin{array}{lcr}
  \ 0 &  -D_{-X}\otimes I_r \\
    \  D_{-X}\otimes I_r  & 0
\end{array}\right];~
p=\left[\begin{array}{lcr}
  \ p &  0 \\
    \ 0  & p
\end{array}\right];~$$
$$
\sigma=i\left[\begin{array}{lcr}
  \ 0 &  I_r \\
    \  I_r  & 0
\end{array}\right];~
e^{-X}=\left[\begin{array}{lcr}
  \  e^{-X} &  0 \\
    \ 0  & e^{-X}
\end{array}\right],$$
be operators from $H\otimes {\bf C^r}\oplus H\otimes {\bf C^r}$ to
itself, then
$${\bf D_{-X}}\sigma=-\sigma {\bf D_{-X}};~~\sigma p=p\sigma.$$
Moreover ${\bf D_{-X}}e^{t{\bf D_{-X}}^2}$ and $e^{t{\bf D_{-X}}^2}~ (t>0)$ are
traceclass. For $u\in [0,1]$, let
$$D_{-X,u}=(1-u)D_{-X}+u[pD_{-X}p+(1-p)D_{-X}(1-p)]=D_{-X}+u(2p-1)[D,p],$$
then
$${\bf
D}_{-X,u}=\left[\begin{array}{lcr}
  \  0 &  -D_{-X,u} \\
    \ D_{-X,u } & 0
\end{array}\right]
={\bf D_{-X}}+u(2p-1)[{\bf D_{-X}},p].$$ \noindent Consider a family of
Dirac Operators parameterized by
$(u,s,t)$, which is given by
$$\widetilde{{\bf D_{-X}}}=t^{\frac{1}{2}}{\bf
D}_{-X,u}+s\sigma(p-\frac{1}{2}).$$ Let $A=d+\widetilde{{\bf D_{-X}}}$ be a
superconnection on the trivial infinite dimensional superbundle with
base $[0,1]\times {\bf R}$ and fibre $H\otimes {\bf
C^r}\oplus H\otimes {\bf C^r}.$ Then we have
$$B_{X,u,s,t}:=(d+\widetilde{{\bf D_{-X}}})^2=t{\bf
D}_{-X,u}^2-s^2/4-(1-u)t^{\frac{1}{2}}s\sigma[{\bf D},p]+ds\sigma
(p-\frac{1}{2})+t^{\frac{1}{2}}du(2p-1)[{\bf
D},p].\eqno(3.40)$$
Consider the differential form $\int_\varepsilon^{+\infty}\frac{1}{2\sqrt{t}}{\rm Str}[\psi_t e^{-tX}{\bf
D}_{-X,u}e^{B_{X,u,s,t}}]_J$ on $[0,1]\times {\bf R}$. By Lemma 9.15 in [BGV] as well as that $D_u$ is inverse, we have
$$d\int_\varepsilon^{+\infty}\frac{1}{2\sqrt{t}}{\rm Str}[\psi_t e^{-tX}{\bf
D}_{-X,u}e^{B_{X,u,s,t}}]_J$$
$$=-\int_\varepsilon^{+\infty}\frac{\partial}{\partial t}{\rm Str}[\psi_t e^{-tX}e^{B_{X,u,s,t}}]_J
={\rm Str}[e^{-X}e^{B_{X/\varepsilon,u,s,\varepsilon}}]_J.\eqno(3.41)$$
 \noindent Let $\Gamma_u=\{u\}\times {\bf R}\subset [0,1]\times {\bf
 R}$ be a contour oriented in the direction of increasing $s$
 and $\gamma_s=[0,1]\times \{s\}$
 be a contour oriented in the direction of increasing $u$ . By the
 Stokes theorem, then
 $$\int_{[0,1]\times {\bf
 R}}d\int_\varepsilon^{+\infty}\frac{1}{2\sqrt{t}}{\rm Str}[\psi_t e^{-tX}{\bf
D}_{-X,u}e^{B_{X,u,s,t}}]_J$$
$$=\left(\int_{\Gamma_1}-\int_{\Gamma_0}
-\int_{\gamma_{+\infty}}+\int_{\gamma_{-\infty}}\right)
\left[\int_\varepsilon^{+\infty}\frac{1}{2\sqrt{t}}{\rm Str}[\psi_t e^{-tX}{\bf
D}_{-X,u}e^{B_{X,u,s,t}}]_J\right].\eqno(3.42)$$
We have for some constant $C>0$ that,
$$\int_{\gamma_s}
\left[\int_\varepsilon^{+\infty}\frac{1}{2\sqrt{t}}{\rm Str}[\psi_t
e^{-tX}{\bf D}_{-X,u}e^{B_{X,u,s,t}}]_J\right]\sim
O(e^{-cs^2}).\eqno(3.43)$$
 As shown in [Ge2] or [Wa1], it can be

$$\int_{\Gamma_0}\left[\int_\varepsilon^{+\infty}\frac{1}{2\sqrt{t}}{\rm Str}[\psi_t e^{-tX}{\bf
D}_{-X,u}e^{B_{X,u,s,t}}]_J\right]$$
$$=
-4\sqrt{-1}\pi[\langle \eta_X(D)_J,{\rm
Ch}(p)\rangle-\frac{1}{2}\langle \eta_X(D)_J,{\rm rk}(p){\rm
Ch}_*(1)\rangle].\eqno(3.44)$$
$$\int_{\Gamma_1}\left[\int_\varepsilon^{+\infty}\frac{1}{2\sqrt{t}}{\rm Str}[\psi_t e^{-tX}{\bf
D}_{-X,u}e^{B_{X,u,s,t}}]_J\right]$$
$$=
-2\sqrt{-1}\pi\left[ \eta_X(D_p)_J-\langle \eta_X(D)_J,{\rm rk}(p){\rm
Ch}_*(1)\rangle\right.$$
$$\left.+\frac{1}{2}\int_0^1{\rm Tr}[\varepsilon^{1/2}(2p-1)dpe^{-X}e^{-\varepsilon D_ {-\frac{X}{\varepsilon},u}^2}]_Jdu\right].\eqno(3.45)$$
By (3.41)-(3.45), we get\\

\noindent {\bf Theorem 3.9} {\it Assume $D$ is inverse and $||dp||<\lambda$ where $\lambda$ is the smallest eigenvalue of $|D|$ and $X$ is small , we have}
 $$\frac{1}{2}\eta_X(p(D\otimes
 I_r)p)_J=\langle\eta_X(D)_J,{\rm Ch}(p)\rangle$$
 $$+\pi\sqrt{-1}\int_0^1{\rm Tr}[\varepsilon^{1/2}(2p-1)dpe^{-X}e^{-\varepsilon D_ {-\frac{X}{\varepsilon},u}^2}]_Jdu
 -\frac{1}{4\sqrt{-1}\pi}\int_{[0,1]\times {\bf
 R}}{\rm
 Str}[e^{-X}e^{B_{X/\varepsilon,u,s,\varepsilon}}]_J.\eqno(3.46)$$

\section{A family infinitesimal equivariant index formula}

 \quad In this section, we give a proof of a family infinitesimal equivariant index formula.
 Let $M$ be a $n+{q}$ dimensional compact connected manifold
 and $B$ be a ${q}$ dimensional compact connected
 manifold. Assume that $\pi :M\rightarrow B$ is a submersion of
 $M$ onto $B$, which defines a fibration of $M$ with fibre $Z$. For
 $y\in B$, $\pi^{-1}(y)$ is then a submanifold $Z_y$ of $M$. Let $TZ$
 denote the $n$-dimensional vector bundle on $M$ whose fibre $T_xZ$
 is the tangent space at $x$ to the fibre $Z_{\pi x}$. Assume
 that $M$ and $B$ are oriented. Taking the orthogonal bundle of $TZ$
 in $TM$ with respect to any Riemannian metric, determines a smooth
 horizontal subbundle $T^HM$, i.e. $TM=T^HM\oplus TZ$. A vector field
 $Y\in TB$ will be identified with its horizontal lift $Y\in
 T^HM$, moreover $T^H_xM$ is isomorphic to $T_{\pi(x)}B$ via
 $\pi_*$. Recall that $B$ is Riemannian, so we can lift the
 Euclidean scalar product $g_B$ of $TB$ to $T^HM$.
And we assume that $TZ$ is endowed with a scalar product $g_Z$. Thus
we can introduce a new scalar product $g_B\oplus g_Z$ in $TM$.
Denote by $\nabla^L$ the Levi-Civita connection on $TM$ with respect
to this metric. Let $\nabla^B$ denote the Levi-Civita connection on
$TB$ and still denote by $\nabla^B$ the pullback connection on
$T^HM$. Let $\nabla^Z=P_Z(\nabla^L)$, where $P_Z$ denotes the
projection to $TZ$. Let $\nabla^{\oplus}=\nabla^B\oplus \nabla^Z$
and $\omega=\nabla^L-\nabla^{\oplus}$ and $T$ be the torsion tensor
of $\nabla^{\oplus}$. Let $SO(TZ)$ be the $SO(n)$ bundle of oriented
orthonormal frames in $TZ$. Now we assume that the bundle $TZ$ is spin.
Let $S(TZ)$ be the associated spinors bundle and $\nabla^Z$ can be
lifted to give a connection on $S(TZ)$. Let $D$ be the tangent Dirac
operator.\\
\indent Let $G$ be a compact Lie group which acts fiberwise on $M$.
We will consider that $G$ acts as identity on $B$. Without loss of
generality assume $G$ acts on $(TZ, h_{TZ})$
isometrically. Also assume that the action of G lifts to $S(TZ)$
and the $G$-action commutes with $D$. Let $E$ be the vector
bundle $\pi^*(\wedge T^*B)\otimes S(TZ)$. This bundle carries a
natural action $m_0$ of the degenerate Clifford module $C_0(M)$. The
Clifford action of a horizontal cotangent vector
$\alpha\in\Gamma(M,T_H^*M)$ is given by exterior multiplication
$m_0(\alpha)=\varepsilon(\alpha)$ acting on the first factor
$\bigwedge T_H^*M$ in $E$, while the Clifford action of a vertical
cotangent vector simply equals its Clifford action on $S(TZ)$.
Define the connection for $X\in {\mathfrak{g}}$ whose Killing vector
field is in $TZ$,
$$\nabla^{E,-X,\oplus}:=\pi^*\nabla^B\otimes 1+1\otimes \nabla^{S,-X},\eqno(4.1)$$
$$\omega(Y)(U,V):=g(\nabla^L_YU,V)-g(\nabla^{\oplus}_YU,V),\eqno(4.2)$$
$$\nabla^{E,-X,0}_Y:=\nabla^{E,-X,\oplus}_Y+\frac{1}{2}m_0(\omega(Y)),\eqno(4.3)$$
for $Y,U,V\in TM$. Then the Bismut Connection acting on $\Gamma(M,\pi^*\wedge(T^*B)\otimes S(TZ))$ is defined by
$${\mathcal{B}}^{-X}=\sum_{i=1}^nc(e_i^*)\nabla_{e_i}^{E,-X,0}+\sum_{j=1}^{{q}}c(f_j^*)\nabla_{f_j}^{E,-X,0},\eqno(4.4)$$
where $e_1,\cdots,e_n$ and $f_1,\cdots,f_q$ are orthonormal basis of $TZ$ and $TB$ respectively.
Define the family Bismut Laplacain as follows:
$$H_X=({\mathcal{B}}^{-X})^2+{\mathfrak{L}}^E_X.\eqno(4.5)$$
Let $\triangle^{Z,X}$ be the Laplacian on $\pi^*(\wedge T^*B)\otimes S(TZ)$ associated with $\nabla^{E,X,0}$. Similarly to Proposition
8.12 and Theorem 10.17 in [BGV], we have\\

\noindent{\bf Proposition 4.1}{\it ~ The following identity holds}
$$H_X=\triangle^{G,X}+\frac{1}{4}r_M.\eqno(4.6)$$\\

In this section, we establish an index theorem in the untwisted case and it is easy to extend it to the twisted case.
Then by Proposition 10.15 in [BGV],
$${\mathcal{B}}^{-X}={\mathcal{B}}+\frac{1}{4}c(X)=D_{-X}+\textrm{A}_{[+]},\eqno(4.7)$$ where $\textrm{A}_{[+]}$ is an
operator with coefficients in $\Omega_{\geq 1}(B)$ and ${\mathcal{B}}$ is the Bismut superconnection. And $D_{-X}=D+\frac{1}{4}c(X)$. Let
$ H_X=D_{-X}^2+F_{[+]},$ where $\textrm{F}_{[+]}$ is an
operator with coefficients in $\Omega_{\geq 1}(B)$. We define the operator $e^{-t{H_X}}$ which is given by
$$e^{-t{H_X}}=e^{-t(D_{-X}^2+L_X)}+\sum_{k>0}(-t)^kI_k,\eqno(4.8)$$
where
$$I_k=\int_{\triangle_k}e^{-\sigma_0t(D_{-X}^2+L_X)}F_{[+]}e^{-\sigma_1t(D_{-X}^2+L_X)}F_{[+]}$$
$$\cdots e^{-\sigma_{k-1}t(D_{-X}^2+L_X)}F_{[+]}e^{-\sigma_kt(D_{-X}^2+L_X)}
d\sigma,\eqno(4.9)$$
and the sum is finite. By Theorem 2.1 in [LM], similarly to Proposition 8.11 in [BGV], we get\\

\noindent{\bf Proposition 4.2} {\it We have in the cohomology class of $B$,}
$${\rm Ch}({\rm ind}_G(e^{-X},D))={\rm Str}(\phi_te^{-tH_{\frac{X}{t}}}),\eqno(4.10)$$
{\it which does not depend on $t$, and $\phi_t(dy_j)=\frac{1}{\sqrt{t}}dy_j$.}\\

We define the operator
$$Q:=(H_X+\frac{\partial}{\partial t})^{-1}=(D_{-X}^2+L_X+\frac{\partial}{\partial t})^{-1}$$
$$+\sum_{k>0}(-1)^k(D_{-X}^2+L_X+\frac{\partial}{\partial t})^{-1}
[F_{[+]}(D_{-X}^2+L_X+\frac{\partial}{\partial t})^{-1}]^k,\eqno(4.11)$$
where $(D_{-X}^2+L_X+\frac{\partial}{\partial t})^{-1}$ is the Volterra
inverse of $D_{-X}^2+L_X+\frac{\partial}{\partial t}$ as shown in Section 2.
 We can define Volterra symbols with coefficients in ${\mathbb{C}}[{\mathfrak{g}}^*]\otimes \wedge T_z^*B$
  and Volterra pseudodifferential operators with coefficients in ${\mathbb{C}}[{\mathfrak{g}}^*]\otimes \wedge T_z^*B$.
   Write the space of Volterra pseudodifferential operators with coefficients in ${\mathbb{C}}[{\mathfrak{g}}^*]\otimes \wedge T_z^*B$ by
  $\Psi_V^*({\mathbb{R}}^n\times {\mathbb{R}}, S(TM)\otimes {\mathbb{C}}[{\mathfrak{g}}^*]\otimes \wedge T_z^*B).$
We define the Getzler order $O_G(dy^j)=1.$ Let $Q\in
\Psi_V^*({\mathbb{R}}^n\times {\mathbb{R}}, S(TM)\otimes
{\mathbb{C}}[{\mathfrak{g}}^*]\otimes \wedge T_z^*B)$ have symbol
$$q(x,X,\xi,\tau)\sim
 \sum_{k\leq m'}\sum_{l_0=0}^{2^{{\rm dim}B}}q_{k,l_0}(x,X,\xi,\tau)\omega_{[l_0]},\eqno(4.12)$$
 where $q_{k,l_0}(x,X,\xi,\tau)$ is an order $k$ symbol and $\omega_{[l_0]}$ is a $l_0$-degree differential form on $B$. Then using Taylor expansions at $x = 0$ and at $X=0$, it gives that
$$\sigma[q(x,X,\xi,\tau)]\sim\sum_{j,k,\alpha,\beta}\sum_{l_0}^{2^{{\rm dim}B}}\frac{x^\alpha}{\alpha!}\frac{X^\beta}{\beta!}
\sigma[\partial_x^\alpha \partial_X^\beta
q_{k,l_0}(0,0,\xi,\tau)]^{(j)}\omega_{[l_0]}.\eqno(4.13)$$ The
symbol $\frac{x^\alpha}{\alpha!}\frac{X^\beta}{\beta!}
\sigma[\partial_x^\alpha \partial_X^\beta
q_{k,l_0}(0,0,\xi,\tau)]^{(j)}\omega_{[l_0]}$ is the Getzler
homogeneous of $k+j+l_0-|\alpha|+2|\beta|$. Similarly to the
definition 2.7, we can define
 the $J$-truncated symbol of $q$ denoted by
 $\sigma[q(x,X,\xi,\tau)]_J.$ Also, we may define
  the truncated model operator of $Q$. Similarly to Lemma 2.9, we have\\

\noindent {\bf Lemma 4.3} {\it Let
$Q\in\Psi_V^*({\mathbb{R}}^n\times {\mathbb{R}}, S(TM)\otimes
{\mathbb{C}}[{\mathfrak{g}}^*]\otimes \wedge T_z^*B),$ and $Q_J$
has the Getzler order $m$ and model operator $Q_{(m),J}$. Then as
$t\rightarrow 0^+$ we have:}
$${\rm 1)}~~ \sigma[\phi_tK_{Q_J}(0,0,\frac{X}{t},t)]^{(j)}=\omega O(t^{\frac{j-n-m-2}{2}})+O(t^{\frac{j-n-m-1}{2}}),~~ {\rm where} ~\omega\in\Omega^{{\rm odd}}(T^*B),
$$
$~~{\rm if~} m-j~~{\rm ~is~ odd};$
$${\rm 2)}~~\sigma[\phi_tK_{Q_J}(0,0,\frac{X}{t},t)]^{(j)}=t^{\frac{j-n-m-2}{2}}K_{Q_{(m),J}}(0,0,X,1)^{(j)}+O(t^{\frac{j-n-m}{2}}),~~{\rm if~} m-j~~{\rm ~is~ even},$$
{\it where $[K_{Q_J}(0,0,\frac{X}{t},t)]^{(j)}$ denotes taking the
$j$ degree form component in $\wedge^*T(M_z)$. In particular, when $m=-2$ and $j=n$
is even, we get}
$$\sigma[\phi_tK_{Q_J}(0,0,\frac{X}{t},t)]^{(n)}=K_{Q_{(-2),J}}(0,0,X,1)^{(n)}+O(t).\eqno(4.14)$$\\

By Lemma 4.3, similarly to the proof of Proposition 2.11 and
Theorem 2.12, we have\\

\noindent{\bf Theorem 4.4} {\it We have in the cohomology class of
$B$,}
$${\rm Ch}({\rm ind}_G(e^{-X},D))={(2\pi\sqrt{-1})}^{-n/2}\int_{M/B}\widehat{A}(F^Z_{\mathfrak{g}}(X))d{\rm
Vol}_{M/B}.\eqno(4.15)$$\\

\indent In the following, we define infinitesimal equivariant eta
forms. Now assume that ${\rm dim}M$ and ${\rm dim}Z$ are odd. Let ${\rm Tr}^{\rm even}$ denote taking trace on
the coefficients of even forms on $B$. Let $T$ be the torsion tensor
of $\nabla^{\oplus}$ and $c(T)=\sum_{1\leq \alpha<\beta \leq
q}dy_\alpha dy_\beta c(T(\frac{\partial}{\partial
y_\alpha},\frac{\partial}{\partial y_\beta})).$ Then the {\bf
infinitesimal equivariant eta form} is defined by
$$\widehat{\eta}_X=\int_0^{\infty}\frac{1}{\sqrt{\pi t}}{\rm Tr}^{\rm
even}[\phi_t(D-\frac{1}{4}c(\frac{X}{t})+\frac{c(T)}{4})e^{-tH_{\frac{X}{t}}}]dt.\eqno(4.16)$$
When $n$ is even, we define the infinitesimal equivariant eta form
by
$$\widehat{\eta}_X=\int_0^{\infty}\frac{1}{\sqrt{\pi t}}{\rm Str}[\phi_t(D-\frac{1}{4}c(\frac{X}{t})+\frac{c(T)}{4})e^{-tH_{\frac{X}{t}}}]dt.\eqno(4.17)$$
 Let $e_1(x),\cdots e_n(x)$ denote the orthonormal frame of $TZ$.
   If $A(Y)$ is any $0$ order operator depending linearly on
   $Y\in TZ$, we define the operator $(\nabla_{e_i}+A(e_i))^2$ as
   follows
   $$(\nabla_{e_i}+A(e_i))^2=\sum_1^n(\nabla_{e_i(x)}+A(e_i(x)))^2-\nabla_{\sum_j\nabla_{e_j}e_j}
   -A(\sum_j\nabla_{e_j}e_j).\eqno(4.18)$$
We introduce an auxiliary form $dt$. Let
$$\widetilde{H_X}=H_X-dt(D_X+\frac{c(T)}{4});~h(x)=1+\frac{1}{2}dt\sum_{j=1}^nx_ic(e_i),\eqno(4.19)$$
then we have
$$\widetilde{H_X}=-(\nabla^{G,X}_{e_i}+\frac{1}{2}<\omega(e_i)e_j,
f_\alpha>e_jdy_\alpha+\frac{1}{4}<\omega(e_i)f_\alpha,f_\beta>dy_\alpha
dy_\beta-\frac{1}{2}c(e_i)dt)^2+\frac{r_Z}{4}.\eqno(4.20)$$ and
$$(h\rho)\widetilde{H_X}(h\rho)^{-1}=\rho H_X\rho^{-1}+dtu,~~ O_G(u)\leq 0 .\eqno(4.21)$$
By the trick in Lemma 3.3, we get\\

\noindent{\bf Lemma 4.5} {\it When $t\rightarrow 0^+$, we have}
$${\rm Tr}^{\rm
even}[\phi_t(D-\frac{1}{4}c(\frac{X}{t})+\frac{c(T)}{4})e^{-tH_{\frac{X}{t}}}]\sim
O(t^{1/2}).\eqno(4.22)$$\\

\noindent{\bf Remark.} We also prove Lemma 4.5 by using the method in [BGV, p. 347].\\

\indent We introduce the following notations as those in Lemma 3.4,
$${\mathcal{M}}_{(j)}:={\mathcal{M}}\cap \oplus_{k+l\geq j}
A^k(B,{\rm
End}(\Gamma(S(TM_z)))\otimes{\mathbb{C}}[{\mathfrak{g^*}}]_{J,l};$$
$$~~{\mathcal{N}}_{(j)}:={\mathcal{N}}\cap \oplus_{k+l\geq j}
A^k(B,{\rm
End}(\Gamma(S(TM_z)))\otimes{\mathbb{C}}[{\mathfrak{g^*}}]_{J,l}.$$
${\mathcal{M}}_{(j)}$ is the algebra generated by differential operators and smoothing operators acting on $\Gamma(S(TM_z))$ with coefficients in
$\oplus_{k+l\geq j}
A^k(B)\otimes{\mathbb{C}}[{\mathfrak{g^*}}]_{J,l}$
and ${\mathcal{N}}_{(j)}$ denotes the algebra generated by smoothing operators acting on $\Gamma(S(TM_z))$ with coefficients in $\oplus_{k+l\geq j}
A^k(B)\otimes{\mathbb{C}}[{\mathfrak{g^*}}]_{J,l}.$
Replace
${\mathcal{M}}_{j}$ and ${\mathcal{N}}_{j}$ in Lemma 2.34 in
[Go] by ${\mathcal{M}}_{(j)}$ and ${\mathcal{N}}_{(j)}$, then we have\\

\noindent{\bf Lemma 4.6} {\it We assume that the kernel of $D$ is a complex vector bundle. When $t\rightarrow +\infty$, we have}
$${\rm Tr}^{\rm
even}[\phi_t(D-\frac{1}{4}c(\frac{X}{t})+\frac{c(T)}{4})e^{-tH_{\frac{X}{t}}}]\sim
O(t^{-1}).\eqno(4.23)$$\\

By Lemma 4.5 and Lemma 4.6, we get that infinitesimal
equivariant eta forms are well-defined. We recall the definition of
equivariant eta forms in [Wa2],
$$\widehat{\eta}(e^{-X})=\int_0^{\infty}\frac{1}{\sqrt{\pi t}}{\rm Tr}^{\rm
even}[\phi_te^{-X}(D+\frac{c(T)}{4})e^{-t{\mathcal{B}}^2}]dt.\eqno(4.24)$$
In the last, we announce a comparison formula between infinitesimal
equivariant eta forms and equivariant eta forms and its proof will
appear elsewhere. Let $d_{rX}\theta_{rX}=d\theta_{rX}-||rX||$ for $r>0$. The vector field $X_M$ is called geodesic if $\nabla_{X_M}X_M=0$.\\

\noindent{\bf Theorem 4.7} {\it If the
Killing field $X_M$ is geodesic and has no zeros on $M$, then for $X\in {\mathfrak{g}}$ and
small $r\neq 0$ and each $K>0$, we have up to an exact form}
$$\widehat{\eta}_{rX}=\widehat{\eta}(e^{-rX})+\int_{M/B}2(2\pi
i)^{-\frac{n+1}{2}}\frac{\theta_{rX}}{d_{rX}\theta_{rX}}\widehat{A}_{rX}(M/B)+O(r^K).\eqno(4.25)$$\\

 \noindent {\bf Acknowledgements.} This work
was supported by NSFC No.11271062 and NCET-13-0721. The author would like to thank Profs.
Weiping Zhang and Huitao Feng for introducing the index theory to
him. He would like to thank Prof. S. Goette for helpful discussions. The author is indebted to the referee for his (her) careful reading and helpful comments.\\

\noindent{\large \bf References}\\

\noindent[Az]F. Azmi, The equivariant Dirac cyclic cocycle, Rocky
Mountain J. Math. 30 (2000), 1171-1206.

\noindent[BGS]R. Beals, P. Greiner and N. Stanton, The heat equation on
a CR manifold, J. Differential Geom. 20 (1984), 343-387.

\noindent[BGV]N. Berline, E. Getzler and M. Vergne, {\it Heat kernels
and Dirac operators,} Springer-Verlag, Berlin, 1992.

\noindent[BV1]N. Berline and M. Vergne, A computation of the
equivariant index of the Dirac operators, Bull. Soc. Math. France
113 (1985), 305-345.

 \noindent[BV2]N. Berline and M. Vergne, The
equivariant index and Kirillov character formula, Amer. J. Math.
107 (1985), 1159-1190.

\noindent [Bi]J. M. Bismut, The infinitesimal Lefschetz formulas: a
heat equation proof, J. Func. Anal. 62 (1985), 435-457.

\noindent[BlF]J. Block and J. Fox, Asymptotic pseudodifferential operators and index theory, Contemp. Math., 105 (1990), 1-32.

\noindent [CH]S. Chern and X. Hu, Equivariant Chern character
for the invariant Dirac operators, Michigan Math. J. 44 (1997),
451-473.

\noindent[Co]A. Connes, Entire cyclic cohomology of Banach algebras and characters of
$\theta$-summable Fredholm module, K-Theory 1 (1988), 519-548.

\noindent[CM1]A. Connes and H. Moscovici, Cyclic cohomology, the Novikov conjecture and hyperbolic groups,
Topology 29 (1990), 345-388.

\noindent [CM2]A. Connes and H. Moscovici, Transgression and
Chern character of finite dimensional K-cycles,
 Commun. Math. Phys. 155 (1993), 103-122.

\noindent[Do1]H. Donnelly, Eta invariants for $G$-spaces, Indiana
Univ. Math. J. 27 (1978), 889-918.

\noindent[Fa]H. Fang, Equivariant spectral flow and a Lefschetz
theorem on odd dimensional spin manifolds, Pacific J. Math.
220 (2005), 299-312.

\noindent[Fe]H. Feng, A note on the noncommutative Chern character
(in Chinese), Acta Math. Sinica 46 (2003), 57-64.

\noindent[Go]S. Goette, Equivariant eta invariants and eta forms, J.
reine angew Math. 526 (2000), 181-236.

\noindent[Ge1]E. Getzler, The odd Chern character in cyclic homology
and spectral flow, Topology 32 (1993), 489-507.

\noindent [Ge2]E. Getzler, Cyclic homology and the
Atiyah-Patodi-Singer index theorem, Contemp. Math. 148 (1993), 19-45.

\noindent [GS]E. Getzler and A. Szenes, On the Chern character of
theta-summable Fredholm modules, J. Func. Anal. 84 (1989), 343-357.

\noindent[Gr]P. Greiner, An asymptotic expansion for the heat equation, Arch. Rational Mech. Anal. 41 (1971), 163-218.

\noindent[JLO]A. Jaffe, A. Lesniewski and K. Osterwalder, Quantum K-theory: The Chern
character, Comm. Math. Phys. 118 (1988), 1-14.

\noindent [KL]S. Klimek and A. Lesniewski, Chern character in
equivariant entire cyclic cohomology, K-Theory 4 (1991), 219-226.

\noindent [LYZ]J. D. Lafferty, Y. L. Yu and W. P. Zhang, A
direct geometric proof of Lefschetz fixed point formulas, Trans.
AMS. 329 (1992), 571-583.

\noindent [LM]K. Liu; X. Ma, On family rigidity theorems, I. Duke
Math. J. 102 (2000), no. 3, 451-474.

 \noindent[Po]R. Ponge, A new short proof of the local index formula and some of its applications,
Comm. Math. Phys. 241 (2003), 215-234.

 \noindent[PW]R. Ponge and H. Wang, Noncommutative geometry, conformal geometry, and the local
equivariant index theorem, arXiv:1210.2032.

\noindent[Wa1]Y. Wang, The equivariant noncommutative
Atiyah-Patodi-Singer index theorem, K-Theory, 37 (2006), 213-233.

\noindent[Wa2]Y. Wang, The Greiner's approach of heat kernel
asymptotics, equivariant family JLO characters and equivariant eta
forms, arXiv:1304.7354.

 \noindent[Wu]F. Wu, The Chern-Connes character for the Dirac
operators on manifolds with boundary, K-Theory 7 (1993), 145-174.

\noindent[Zh]W. Zhang, A note on equivariant eta invariants, Proc. AMS 108 (1990), 1121-1129.\\

 \indent{\it School of Mathematics and Statistics, Northeast Normal University, Changchun Jilin, 130024, China }\\
 \indent E-mail: {\it wangy581@nenu.edu.cn}\\

\end{document}